\documentclass[11pt,psamsfonts]{amsart}
\usepackage{palatino}
\usepackage{euler,color}
\usepackage[normalem]{ulem}

\usepackage[all]{xy}
\usepackage{amssymb,amscd,amsmath,amsthm,dsfont}

\chardef\bslash=`\\ 

\newtheorem{theorem}{Theorem}[section]
\newtheorem{cor}[theorem]{Corollary}
\newtheorem{lemma}[theorem]{Lemma}
\newtheorem{proposition}[theorem]{Proposition}
\newtheorem{prop}[theorem]{Proposition}

\theoremstyle{remark}
\newtheorem{remark}[theorem]{Remark}

\newtheorem{example}[theorem]{Example}

\theoremstyle{definition}
\newtheorem{definition}[theorem]{Definition}

\numberwithin{equation}{section}

\newcommand{\proref}[1]{Proposition~\ref{#1}}

\newcommand{\clsp}{\overline{\operatorname{span}}}
\newcommand{\lsp}{\operatorname{span}}
 \newcommand{\id}{\operatorname{id}}
\newcommand{\inver}{\operatorname{inv}}
 \newcommand{\Aut}{\operatorname{Aut}}
\newcommand{\End}{\operatorname{End}}
 
\newcommand{\mult}{\operatorname{mult}}
\newcommand{\add}{\operatorname{add}}
\newcommand{\bdry}{\operatorname{CL}}
\newcommand{\rt}{{\operatorname{rt}}}

\newcommand{\nx}{\mathbb N^{\times}}

\newcommand{\N}{\mathbb N}

\newcommand{\Z}{\mathbb Z}
\newcommand{\Q}{\mathbb Q} 
\newcommand{\qx}{\mathbb Q^\times_+}
\newcommand{\qp}{{\mathbb Q_+}}

\newcommand{\C}{\mathbb C}
\newcommand{\R}{\mathbb R}
\newcommand{\T}{\mathbb T}
\newcommand{\CC}{\mathcal C}
\newcommand{\TT}{\mathcal T}
\newcommand{\F}{\mathcal F}

\newcommand{\NT}{\mathcal{NT}}
\newcommand{\NO}{\mathcal{NO}}

\newcommand{\primes}{\mathcal P}
\newcommand{\PP}{\mathcal P}

\def\lcm{\operatorname{lcm}}
\newcommand{\Prim}{\operatorname{Prim}}
\newcommand\qps{\Q^*_+}

\newcommand{\bqr}{B_{\Q}(\R)}
\def\op{\operatorname{op}}

\newcommand{\nxxn}{{\mathbb N^\times}\!\! \ltimes \mathbb N}

\newcommand{\qxxq}{{\mathbb Q_+^\times}\!\! \ltimes \mathbb Q}
\newcommand{\nxnx}{{\mathbb N \rtimes \mathbb N^\times}}
\newcommand{\qxqx}{{\mathbb Q \rtimes \mathbb Q^\times_+}}

\title[A right Toeplitz algebra]{Boundary quotients of the right Toeplitz algebra of\\
the affine semigroup over the natural numbers}
\date{4 March 2021, with minor changes 13 August 2021}

\thanks{This research has been supported by the Marsden Fund of the Royal Society of New Zealand and the Natural Sciences and Engineering Research Council of Canada Discovery Grant RGPIN-2017-04052.}

\author[Astrid~an~Huef]{Astrid an Huef}
\address{Astrid an Huef and Iain Raeburn\\ 
School of Mathematics and Statistics\\Victoria University of Wellington\\PO Box 600\\Wellington 6140\\New Zealand}
\email{astrid.anhuef@vuw.ac.nz, iain.raeburn@vuw.ac.nz}
\author[Marcelo~Laca]{Marcelo Laca}
\address{Marcelo Laca\\ Department of Mathematics and Statistics\\University of Victoria\\PO~Box~1700\\Victoria\\BC V8W 2Y2\\Canada}
\email{laca@uvic.ca}
\author[Iain~Raeburn]{Iain Raeburn}


\begin{document}

\begin{abstract} We study the Toeplitz $C^*$-algebra   generated by the  right-regular representation of the semigroup $\nxnx$, which we call the right Toeplitz algebra.  We analyse its structure by studying three distinguished quotients. We show that the multiplicative boundary quotient is isomorphic to a crossed product of the Toeplitz algebra of the additive rationals by an action of the multiplicative rationals, and study its ideal structure. The Crisp--Laca boundary quotient is isomorphic to the $C^*$-algebra of the group $\qxxq$.  There is a natural dynamics on the right Toeplitz algebra and all its   KMS states   factor through the additive boundary quotient. We  describe the KMS simplex for inverse temperatures greater than one. 
\end{abstract}
\maketitle

\section*{We dedicate this paper to the memory of Vaughan Jones}
We remember Vaughan with great admiration and fondness. More specifically, in chronological order:
\subsection*{Iain} I met Vaughan at a conference in Kingston, Ontario, shortly after he finished his PhD in 1980; just four years later, we participated in a year-long program at MSRI, and Vaughan was being enthusiastically recruited by Berkeley. I remember a lovely dinner in 1990 with Canadian colleagues, where we toasted Vaughan's Fields medal with an Australian red; we think he would have approved. But it was only when I moved here in 2009 that I learned about all the wonderful things he had been doing for mathematics in New Zealand. We may have had only a little bit of his time, but he used it to do an awful lot of good. It has been a pleasure and an honour to know him. 
\subsection*{Marcelo}
I remember distinctly to this day that in 1984, at a seminar in MSRI, a complete stranger sat next to me and proceeded to fill pages after pages with incomprehensible little scribbles that had numerical labels underneath. I was a beginning graduate student then and found this utterly bewildering, but, Berkeley being Berkeley, I thought maybe he was there just for the coffee and cookies. Until I realized that the scribbles were knots and that sitting next to me was V. F. R. Jones in the very act of classifying them. He was recruited by UC Berkeley the following year so I got to know him better, first as a professor and in due time as a colleague. He was relaxed and intense, supportive and challenging, all at the same time. Very naturally over the years we became friends; his unassuming and easygoing demeanor and his enthusiasm made that unavoidable. His mathematical genius and his personal flair sprinkled with irreverence set him aside as one of a kind. Vaughan also transcended the boundaries of the discipline in many directions. I have had school kids excitedly tell me about something really cool that a certain Mr Jones had done with knots, only to see their jaw drop when I tell them that I knew Mr Jones myself. It was a privilege to have shared with Vaughan what, in retrospect, were far too few precious moments. He will be remembered and he will be missed.
\subsection*{Astrid}
Everyone agrees that Vaughan has done wonderful mathematics and wonderful things for mathematics all over the world. But I really treasure how generous Vaughan was with his time and energy, and in particular how kind and supportive he was of me when I first moved to NZ in 2010. His influence on the NZ mathematics community has been colossal.

\section{Introduction}
We consider the semidirect product group $\qxqx$ and its Toeplitz algebras. In previous work \cite{LR-advmath}, we studied the usual Toeplitz algebra $\TT(\nxnx)$ associated to the left-invariant partial order with positive cone $\nxnx$. The main results there show that there is a natural dynamics on $\TT(\nxnx)$ that admits a rich supply of KMS states and exhibits phase transitions like those of the Bost--Connes system \cite{bos-con}.

The semigroup $\nxnx$ also induces a right-invariant partial order on $\qxqx$, but the corresponding Toeplitz algebra $\TT_{\text rt}(\nxnx)$ is known to be quite different. For example, the left Toeplitz algebra $\TT(\nxnx)$ has an ideal $I$ which contains every other proper ideal, and for which the boundary quotient $\TT(\nxnx)/I$ is  the simple purely infinite $C^*$-algebra $\mathcal{C}_{\Q}$ of Cuntz \cite{cun2}. The right Toeplitz algebra $\TT_{\text rt}(\nxnx)$, on the other hand, has many one-dimensional representations. Curiously, though, in spite of this seemingly huge difference, the two Toeplitz algebras are known to have the same $K$-theory \cite[\S6.4]{CEL}.

Our goal here is to study  the structure of the right Toeplitz algebra $\TT_{\text rt}(\nxnx)$. As suggested in \cite{CEL}, we view it as the left Toeplitz algebra of the opposite semigroup $(\nxnx)^{\op}$, because we can then use the general results about left-invariant quasi-lattice ordered groups from \cite{nica} and \cite{quasilat}. We then avoid a proliferation of $\op$s by realising $(\nxnx)^{\op}$ as the isomorphic semidirect product $\nxxn$. This notation leads us to use similar notation for crossed-product $C^*$-algebras, and we discuss some of the ramifications in \S\ref{cps}.

The Toeplitz algebra $\TT(\nxxn)$ is generated by an isometric representation $V$ of $\nx$ and a single isometry $S$ which generates a representation of $\N$. We use these generators to give a presentation of $\TT(\nxxn)$ (Proposition~\ref{presentation}). Then we focus on three quotients of $\TT(\nxxn)$: the multiplicative boundary quotient $\partial_{\mult}\TT(\nxxn)$ in which the $V$s are unitary, the additive boundary quotient $\partial_{\add}\TT(\nxxn)$ in which $S$ is unitary, and the  Crisp--Laca boundary quotient $\partial\TT(\nxxn)$ in which they all are. 

Our first main result is a structure theorem for $\partial_{\mult}\TT(\nxxn)$. For the proof, we realise  $\partial_{\mult}\TT(\nxxn)$ as a crossed product of the usual Toeplitz algebra $\TT(\qp)$ by an action of the group $\qx$ (Proposition~\ref{idcpToe}), and analyse the structure of this crossed product $C^*$-algebra. The structure of the Toeplitz algebra $\TT(\qp)$ is described by results of Douglas \cite{dou}. We identify the commutator ideal in the Toeplitz algebra as a corner in a crossed product of a commutative $C^*$-algebra, and then combine the actions to get an action of $\qx\ltimes \Q$ on the commutative algebra, which we can study using the well-established (i.e., now relatively old-fashioned!) theory of transformation group algebras. The final result (Theorem~\ref{structuremultbdary}) describes a composition series for $\partial_{\mult}\TT(\nxxn)$ with a large commutative quotient, and two other simple subquotients.

Then, following \cite{LR-advmath}, we  consider a dynamics $\sigma:\R\to \Aut\TT(\nxxn)$ such that $\sigma_t$ fixes $S$ and $\sigma_t(V_a)=a^{it}V_a$, and aim to describe the KMS states of $(\TT(\nxxn),\sigma)$. We have been successful for inverse temperatures $\beta>1$, by applying general results from \cite{diri} about KMS states on semidirect products in \cite[Theorem~7.1]{LR-advmath}; they are parametrised by probability measures on $\T$ (see Theorem~\ref{mainKMSthm} below). At $\beta=1$, the arguments break down: we can see by example (using various combinations of point masses on $\T$) that there are many KMS$_1$ states, but we do not know how to parametrise them. 

In Appendix~\ref{rabbit} we discuss how our results intersect with those of \cite{kaka} and \cite{ABLS}, and, in particular, realise the additive boundary quotient as a Nica--Toeplitz--Pimsner algebra analogous to the ones studied in \cite{kaka}.

\section{Quasi-lattice ordered groups}

Suppose that $G$ is a group and $P$ is a submonoid which generates $G$ and satisfies $P\cap P^{-1}=\{e\}$. Then there is a partial order $\leq_l$ on $G$ such that $g\leq_l h\Longleftrightarrow g^{-1}h\in P$; the subscript in $\leq_l$ is there to remind us that the partial order is invariant under left multiplication by elements of $G$. Following Nica \cite{nica}, the pair $(G,P)$ is \emph{quasi-lattice ordered} if  
\begin{itemize}
\item[(QL)]  for all $n\geq 1$ and all $x_1, \dots, x_n$ in $G$ with a common  upper bound in $P$ also have a least common upper bound in $P$. 
\end{itemize} 
It follows from  \cite[Lemma~7]{CL1} that   (QL) is equivalent to 
\begin{itemize}
\item[(QL1)] every element of $G$ with an upper bound in $P$ has a least upper bound in $P$.
\end{itemize} 
The subset of elements of $G$ which have upper bounds in $P$ (that is, to which (QL1) applies) is
\[
PP^{-1}=\big\{pq^{-1}:p,q\in P\big\}.
\]
An isometric representation of $P$ is a homomorphism $V$ of $P$ into the semigroup of isometries in a $C^*$-algebra, and $V$ is \emph{Nica covariant} if the range projections satisfy $V_pV_p^*V_qV_q^*=V_{p\vee q}V_{p\vee q}^*$ for all $p,q\in P$. A key example of a Nica-covariant representation  is the Toeplitz representation $T:P\to B(\ell^2(P))$, which is characterised in terms of the usual basis $\{\epsilon_p:p\in P\}$ by $T_q\epsilon_p=\epsilon_{qp}$. The Toeplitz algebra $\TT(P)$ is the $C^*$-subalgebra of $B(\ell^2(P))$ generated by $\{T_p:p\in P\}$.

The pair $(G,P)$ also has a right-invariant partial order such that $g\leq_r h\Longleftrightarrow hg^{-1}\in P$. Then there is a homomorphism $W:P\to B(\ell^2(P))$ characterised by
\[
W_q\epsilon_p=\begin{cases}\epsilon_{pq^{-1}}&\text{if $q\leq_r p$ and}\\
0&\text{otherwise.}\end{cases}
\]
Each $W_q$ is a \emph{coisometry}: $W_q^*:\epsilon_p\to \epsilon_{pq}$ is an isometry. We define the \emph{right Toeplitz algebra} $\TT_\rt(P)$ to be the $C^*$-subalgebra $C^*(W_p:p\in P)=C^*(W^*_p:p\in P)$ of $B(\ell^2(P))$. There is a dual notion of right quasi-lattice order, and a corresponding notion of Nica covariant co-representations in which the initial projections satisfy $W_p^*W_pW_q^*W_q=W_{p\vee q}^*W_{p\vee q}^*$ for all $p,q$. The Toeplitz co-representation $W$ is Nica covariant.

We can study the right Toeplitz algebra $\TT_\rt(P)$ by applying  existing theory to the opposite semigroup $(G^{\op},P^{\op})$, which is the set $G^{\op}=\{g^{\op}:g\in G\}$ with $g^{\op}h^{\op}=(hg)^{\op}$. Indeed, the map $g\mapsto g^{\op}$ is an anti-isomorphism which maps $P$ to $P^{\op}$ and satisfies 
\[
g^{\op}\leq_l h^{\op}\Longleftrightarrow (g^{\op})^{-1}h^{\op}\in P\Longleftrightarrow (hg^{-1})^{\op}\Longleftrightarrow hg^{-1}\in P\Longleftrightarrow g\leq_r h.
\]
 Then the unitary $U:\ell^2(P)\to \ell^2(P^{\op})$ characterised by $U\epsilon_p=\epsilon_{p^{\op}}$ satisfies
\[
UW_q^*U^*\epsilon_{p^{\op}}=UW_q^*\epsilon_p=U\epsilon_{pq}=\epsilon_{q^{\op}p^{\op}}=T_{q^{\op}}\epsilon_{p^{\op}}.
\]
Thus conjugation by $U$ gives an anti-isomorphism of $\TT_\rt(P)=C^*(W_q)$ onto $\TT(P^{\op})=C^*(T_q)$, which is the object we study in this paper.

\section{The opposite of the affine semigroup and its Toeplitz algebra}

We use $\N$ for the natural numbers which for us include $0$. We write $\nx$ for the multiplicative semigroup $\{1,2,3,\dots\}$ of $\N$. We write $\Q$ for the rational numbers,  $\qp$ for the additive semigroup of non-negative rational numbers and $\qx$ for the multiplicative group of positive rational numbers.

We recall from \cite{LR-advmath} that the affine group $\qxqx$ consists of pairs $(r,a)\in \Q\times \qx$ with
\[
(r,a)(s,b)=(r+sa,ab)\quad\text{and}\quad (r,a)^{-1}=(-a^{-1}r,a^{-1}).
\]
We will avoid a proliferation of $\op$s by viewing $(\qxqx)^{\op}$ as $\{(a,r)\in \qx\times \Q\}$ with
\[
(a,r)(b,s)=(ab,br+s)\quad\text{and}\quad (a,r)^{-1}=(a^{-1},-a^{-1}r),
\]
which we denote by $\qxxq$. The positive subsemigroup $(\nxnx)^{\op}$ is then
\[
\nxxn:=\{(a,m):a\in \nx,b\in \N\}\quad\text{with}\quad (a,m)(b,n)=(ab, bm+n).
\]
So the right Toeplitz algebra of $(\qxqx,\nxnx)$ is anti-isomorphic to the usual (left) Toeplitz algebra $\TT(\nxxn)$ acting on $\ell^2(\nxxn)$, and we are interested in $\TT(\nxxn)$.

We begin by checking that the established theory applies to $\TT(\nxxn)$. We recall from \cite[Definition~26]{CL1} that $(G,P)$ is \emph{lattice ordered} every element of $G$ has a least upper bound in $P$. By \cite[Lemma~27]{CL1}, $(G,P$ is  lattice ordered if and  only if $G=PP^{-1}$ and $(G,P)$ is quasi-lattice ordered.

\begin{prop}\label{Plo}
The pair $(\qxxq,\nxxn)$ is lattice ordered.
\end{prop}

For comparison,  the pair $(\qxqx, \nxnx)$ is quasi-lattice ordered \cite[Proposition~2.2]{LR-advmath}  but is not lattice ordered. 

\begin{proof}[Proof of Proposition~\ref{Plo}]
We take elements $(a,r)$ and $(b,s)$ of $\qxxq$. Then, by analogy with the discussion at the foot of page~648 in \cite{LR-advmath}, we have:
\begin{align}\label{leqequiv}
(a,r)\leq(b,s)&\Longleftrightarrow (a,r)^{-1}(b,s)\in \nxxn\\
&\Longleftrightarrow (a^{-1},-a^{-1}r)(b,s)\in \nxxn\notag\\
&\Longleftrightarrow (a^{-1}b, s-a^{-1}rb)\in \nxxn\notag\\
&\Longleftrightarrow a^{-1}b\in \nx\text{ and }s-a^{-1}br\in \N.\notag
\end{align}

We now aim to verify  that $(a,r)$ has a least upper bound in $\nxxn$. We write $a=c^{-1}d$ and $r=c^{-1}k$ with $d,c\in \nx$ and $k\in \Z$, and with common denominator $c\in \nx$ as small as possible. Then from \eqref{leqequiv} we have
\begin{align*}
(a,r)\leq (d,n)&\Longleftrightarrow a^{-1}d\in \nx\text{ and }n-a^{-1}dr\in \N\\
&\Longleftrightarrow a^{-1}d=(c^{-1}d)^{-1}d=c\in\nx\quad\text{and}\quad n-a^{-1}dr=n-cr\in\N\\
&\Longleftrightarrow n\geq cr=k\quad\text{because we know that $c\in \nx$.}
\end{align*}
This says, first, that for $n\geq cr=k$, $(d,n)$ is an upper bound for $(a,r)$. 
It also suggests that with $n:=\max(0,k)$ (remember that $k$ is negative if $r$ is), $(d,n)$ is a candidate for the least upper bound. 

To confirm that $(d,n)$ is a \emph{least} upper bound, suppose that $(b,p)\in \nxxn$ and $(a,r)\leq (b,p)$. Then \eqref{leqequiv} implies that
\[
f:=a^{-1}b\in\nx\quad\text{and}\quad l:=p-a^{-1}br\in\N.
\]
But then $f$ satisfies
\[
a=bf^{-1}\in \nx\quad\text{and}\quad r=(p-l)f^{-1}\in \Z.
\]
Thus ``$c$ as small as possible'' implies that there exists $e\in \nx$ such that $f=ce$. Now we have
\[
d^{-1}b=(d^{-1}a)(a^{-1}b)=c^{-1}f=e\in \nx.
\]
Note that $n$ is either $0$ or $k$. Since $p\in \N$,  we trivially have $p-(d^{-1}b)0\in\N$. So we suppose that $n=k$. Then
\[
p-d^{-1}bn=p-d^{-1}bk=p-d^{-1}bcr=p-ecr=p-fr=l\in \N.
\]
Thus we have $(d,n)\leq (b,p)$, as required. 
\end{proof}

Proposition~\ref{Plo} allows us to apply the existing results about left-invariant quasi-lattice ordered groups to $(\qxxq,\nxxn)$. In particular, since $\qxxq$ is a semidirect product of abelian groups, it is amenable, and hence $(\qxxq,\nxxn)$ is amenable as a quasi-lattice ordered group (see \cite[\S4.5]{nica}). Thus we can apply  \cite[\S4.2]{nica} (or \cite[Corollary~3.8]{quasilat}) to $(\qxxq,\nxxn)$. This gives:

\begin{cor}\label{upToeplitz}
The Toeplitz algebra $(\TT(\nxxn), T)$ is universal for Nica-covariant representations of $(\qxxq,\nxxn)$.
\end{cor}

The left-regular representation of $\nxxn$ on $\ell^2(\nxxn)$ is characterised in terms of its action on point masses by
\[
T_{(a,m)} \varepsilon_{(b,n)} = \varepsilon_{(ab,bm+n)}.
\]
We now give a  presentation of the Toeplitz algebra of $\nxxn$.

\begin{proposition}\label{presentation}
Write 
$S:= T_{(1,1)}$ and $V_a:= T_{(a,0)}$ for $a\in \nx$. 
Then
\begin{itemize}
\item[(T0) ] $S^* S = 1 = V_a^* V_a$,
\item[(T1) ] $SV_a = V_a S^a$,
\item[(T2) ] $V_a V_b = V_{ab}$,
\item[(T3) ] $V_a^* V_b = V_b V_a^*$ when $\gcd(a,b) =1$,
\item[(T4) ] $S^*V_a = V_a {S^*}^a$.
\end{itemize}
Let $C^*(s,v)$ be the universal C*-algebra generated by elements $s$ and $\{v_a: a\in \nx\}$ 
satisfying the relations \textnormal{(T0)--(T4)}. Then there is an isomorphism $\pi$ of $C^*(s,v)$ onto the Toeplitz algebra $\TT(\nxxn)$ such that $\pi(s)=S$ and $\pi(v_a)=V_a$.
\end{proposition}

\begin{proof} 
The relation (T0) says that the generators are isometries. Relation (T1) holds because
\[
SV_a=T_{(1,1)}T_{(a,0)}=T_{(1,1)(a,0)}=T_{(a,a)}=T_{(a,0)}T_{(1,a)}=T_{(a,0)}T_{(1,1)}^a=V_aS^a,
\]
and (T3) because $(a,0)(b,0)=(ab,0)$. Relation (T3) follows from Nica-covariance for the pair
$(a,0)$ and $(b,0)$, and we claim that (T4) follows from Nica covariance for $(1,1)$ and $(a,0)$. Indeed, we have $(1,1)\vee(a,0)=(a,a)$, and $T_{(a,a)}=V_aS^a$. Thus Nica covariance gives
\begin{align*}
S^*V_a&=S^*(SS^*V_aV_a^*)V_a=S^*(T_{(1,1)}T_{(1,1)}^*T_{(a,0)}T_{(a,0)}^*)V_a^*=S^*(T_{(a,a)}T_{(a,a)}^*)V_a\\
&=S^*(V_aS^aS^{*a}V_a^*)V_a=S^*(V_aS^a)S^{*a}=S^*(SV_a)S^{*a}\quad\text{using (T1)}\\
&=V_aS^{*a}.
\end{align*}

Now suppose $S$ and $\{V_a:a\in\nx\}$ are elements of a $C^*$-algebra $A$ satisfying (T0)--(T4), and define $L:\nxxn\to A$ by $L_{(a,m)}=V_aS^m$. Then from (T1) we have
\[
L_{(a,m)}L_{(b,n)}=V_aS^mV_bS^n=V_a(V_bS^{bm})S^n=V_{ab}S^{bm+n}=L_{(ab,bm+n)}.
\]
Thus $L$ is a homomorphism of $\nxxn$ into the semigroup of isometries in $A$. We next check that $L$ is Nica covariant. In the following calculations, we write $a=\gcd(a,b)a'$ and similarly for $b$. Then we have on one hand,
\begin{align*}
L_{(a,m)}L_{(a,m)}^*L_{(b,n)}L_{(b,n)}^*&=V_aS^mS^{*m}V_a^*V_bS^nS^{*n}V_b^*=V_aS^mS^{*m}V_{b'}V_{a'}^*S^nS^{*n}V_b^*\\
&=V_{ab'}S^{b'm}S^{*b'm}S^{a'n}S^{*a'n}V_{a'b}^*,
\end{align*}
and on the other
\begin{align*}
L_{(a,m)\vee (b,n)}L_{(a,m)\vee (b,n)}^*&=V_{\lcm(a,b)}S^{\lcm(a,b)(a^{-1}m\vee b^{-1}n)}S^{*\lcm(a,b)(a^{-1}m\vee b^{-1}n)}V^*_{\lcm(a,b)}.
\end{align*}
These last two expressions are the same because (for example) $\lcm(a,b)=\gcd(a,b)a'b'=ab'$, $S^kS^{*k}S^lS^{*l}=S^{k\vee l}S^{*(k\vee l)}$, and $\lcm(a,b)(a^{-1}m\vee b^{-1}n)=b'm\vee a'n$ by left invariance of the partial order. Thus $L$ is Nica covariant.

Since the elements $s$ and $\{v_a:a\in\nx\}$ satisfy the relations (T0)--(T4), the previous paragraph gives us a Nica-covariant representation $L:\nxxn\to C^*(s,v)$ such that $L_{(a,m)}=v_as^m$ for all $(a,m)$. Now Corollary~\ref{upToeplitz} gives a homomorphism $\pi_L:\TT(\nxxn)\to C^*(v,s)$ such that $\pi_L(T_{(a,m)})=v_as^m$. On the other hand, the first paragraph of the proof shows that $S$ and the $V_a$ satisfy the relations in $C^*(v,s)$, and hence the universal property of $C^*(v,s)$ gives a homomorphism of $C^*(s,v)$ into $\TT(\nxxn)$ which is an inverse for $\pi_L$.
\end{proof}

Relations (T0)-(T4) imply that 
\begin{equation}\label{spanning}
\TT(\nxxn)=\clsp\{V_aS^m(V_bS^n)^*: a,b\in\nx, m,n\in\N\}.
\end{equation}

\begin{definition}\label{defn-quotients}
The \emph{multiplicative boundary quotient} $\partial_{\mult}\TT(\nxxn)$ is  the quotient of $\TT(\nxxn)$ obtained by imposing the extra relations 
\begin{itemize}
\item[(T5)] $V_aV_a^*=1$ for all $a\in \nx$.
\end{itemize}
The \emph{additive boundary quotient} $\partial_{\add}\TT(\nxxn)$ is the quotient of $\TT(\nxxn)$ obtained by imposing  the extra relation 
\begin{itemize}
\item[(T6)] $SS^*=1$.
\end{itemize}
The  \emph{Crisp-Laca boundary quotient} $\partial\TT(\nxxn)$ is the quotient of $\TT(\nxxn)$ obtained by imposing the relations (T5) and (T6).
We write $q_{\mult}$, $q_{\add}$ and $q_{\bdry}$ for the quotient maps of $\TT(\nxxn)$ onto $\partial_{\mult}\TT(\nxxn)$, $\partial_{\add}\TT(\nxxn)$ and $\partial\TT(\nxxn)$, respectively. It is possible to verify directly that this is indeed the boundary quotient from \cite{CL2}, but we compute it explicitly below in \proref{propCL=cp}. 
\end{definition}

\begin{remark}  Let $(G,P)$ be a quasi-lattice ordered group. In the literature, for example in  \cite{CL2} and \cite{BaHLR}, the additive, multiplicative and boundary quotients of $\TT(P)$ are defined using the canonical isomorphism of  $\TT(P)$ onto a crossed product $C(X)\rtimes G$ of a partial action of $G$ on a compact space $X$. The quotients correspond to distinguished closed and invariant subsets of $X$. 
\end{remark}
 
\begin{prop}\label{propCL=cp} The Crisp--Laca boundary quotient $\partial\TT(\nxxn)$ is isomorphic to the group $C^*$-algebra $C^*(\qxxq)$. 
\end{prop}
\begin{proof} We will show that $\partial\TT(\nxxn)$ has the universal property of the group $C^*$-algebra $C^*(\qxxq)$. Let $\pi: \partial\TT(\nxxn)\to B(H)$ be a nondegenerate representation on some Hilbert space $H$. Then $\pi=\pi_L$ for some Nica-covariant representation  $L: (\nxxn)\to B(H)$.  The relations (T5) and (T6) imply that $L$ takes values in the group of unitary operators on $H$. Since $(\qxxq,\nxxn)$ is lattice ordered, we have $(\nxxn)(\nxxn)^{-1}=\qxxq$, and $L$ extends to a unitary representation of $\qxxq$ by \cite[Theorem~1.2]{Laca-JLMS}.  Thus every representation of $\partial\TT(\nxxn)$ is induced by a unitary representation of $\qxxq$, as required. Thus $\partial\TT(\nxxn)$ and $C^*(\qxxq)$ are isomorphic. 
\end{proof}

\begin{example}\label{exbdaryrep}
Consider the Hilbert space $\ell^2(\Q_+)$ with orthonormal basis $\{e_x:x\in\qp\}$. Define $Se_x=e_{x+1}$ and $V_ae_x=e_{a^{-1}x}$ for $a\in\nx$. Then we compute 
\[
V_a^*e_x=e_{ax}\quad\text{and}\quad S^*e_x=\begin{cases}0&\text{if $x<1$, and}\\
e_{x-1}&\text{if $x\geq 1$.}\end{cases}.
\]
We think that mental checks suffice for (T1)--(T3) and (T5), but (T4) is more interesting: 
\[
S^*V_ae_x=S^*e_{a^{-1}x}=\begin{cases}0&\text{if $a^{-1}x<1$, and}\\
e_{a^{-1}x-1}&\text{if $a^{-1}x\geq 1$}\end{cases}
\]
is equal to
\[
V_aS^{*a}e_x=\begin{cases}0&\text{if $x<a$, and}\\
e_{a^{-1}(x-a)}=e_{a^{-1}x-1}&\text{if $x\geq a$}\end{cases}
\]
because $a^{-1}x<1\Longleftrightarrow x<a$. Thus we have a representation $\pi_{\mult}$ of $\partial_{\mult}\TT(\nxxn)$ on $\ell^2(\Q_+)$.
\end{example}

\begin{prop}\label{presentaddquot}
Resume the notation of Proposition~\ref{presentation}. Then the images of $S$ and $\{V_a:a\in \nx\}$ in $\partial_{\add}\TT(\nxxn)$ are isometries satisfying:
\begin{itemize}
\item[(A1)] $SV_a=V_aS^a$,
\item[(A2)] $a\mapsto V_a$ is a Nica-covariant representation of $\nx$,
\item[(A3)] $SS^*=1$.
\end{itemize}
The $C^*$-algebra $\partial_{\add}\TT(\nxxn)$ is universal for $C^*$-algebras generated by isometries $s$ and $\{v_a:a\in \nx\}$ satisfying the (lower-case analogues of) these relations.
\end{prop}

\begin{proof}
The relation (T0) in Proposition~\ref{presentation} just says that the elements are isometries, and $(A1)$ follows from (T1). The relations (T2) and (T3) say precisely that $V$ is a Nica-covariant representation, and hence give (A2). The relation (A3) holds because we modded out the ideal generated by $1-SS^*$. 

To see that $\partial_{\add}\TT(\nxxn)$ is universal for these relations, we suppose we have $R,W_a$ satisfying them. Since (A3) implies that $R$ is unitary, we have
\[
RW_a=W_aR^a\Longrightarrow R^*(RW_a)R^{*a}=R^*(W_aR^a)R^{*a}\Longrightarrow W_aR^{*a}=R^*W_a,
\]
and hence the relation (T4) is also satisfied in $C^*(R,W)$. Thus Proposition~\ref{presentation} gives a homomorphism of $\TT(\nxxn)$ into $C^*(R,W)$. Since $R$ is unitary, this homomorphism factors through a homomorphism of $\partial_{\add}\TT(\nxxn)$.
\end{proof}

\begin{example}\label{exDFKcovrep}
We consider the Hilbert space $\ell^2(\nx\times\Z)$ with the usual orthonormal basis 
\[
\{e_{b,m}:b\in \nx, m\in \Z\},
\]
and define $Se_{b,m}= e_{b,b+m}$ and $V_ae_{b,m}=e_{ab,m}$. Then $V$ is a Nica covariant representation of $\nx$, and we have
\begin{align*}
SV_ae_{b,m}&=Se_{ab,m}=e_{ab,ab+m},\quad\text{and}\\
V_aS^ae_{b,m}&=V_ae_{b, ab+m}=e_{ab,ab+m}=SV_ae_{b,m}.
\end{align*}
So (A1) is satisfied. The usual argument shows that $V$ is Nica covariant, and $S$ is unitary because we allow $m$ to be a negative integer. 
\end{example}

\section{Backwards crossed products}\label{cps}
Our choice of the conventions $\nxxn$ and $\qxxq$ for semidirect products, which we use because it allows us to apply results about the usual left-invariant partial order on $(G,P)$, leads us to adopt slightly unorthodox conventions for crossed products of $C^*$-algebras. Suppose that $G$ is a group (in this paper, always discrete) and $\alpha:G\to \Aut A$ is an action of $G$ on a $C^*$-algebra $A$. The \emph{crossed product} is then the $C^*$-algebra generated by a universal covariant representation $(i_A,i_G)$ of the system $(A,G,\alpha)$. As usual, there is up to isomorphism exactly one such crossed product. Here we denote it by $G\ltimes_\alpha A$, suppress the representation $i_A$ of $A$, write $u=i_G$, and view the crossed product as 
\[
G\ltimes_\alpha A:=\clsp\{u_ga:g\in G, a\in A\}.
\]

We now suppose that we have a right action $(q,n)\mapsto n\cdot q$ of one group $Q$ by automorphisms of another group $N$, and $Q\ltimes N$ is the semidirect product with $(p,m)(q,n):=(pq, (m\cdot q)n)$. The unitary representations of $Q\ltimes N$ are given by pairs of representations $U:Q\to U(H)$ and $V:N\to U(H)$ satisfying $U_mV_p=V_pU_{m\cdot p}$. The following lemma for iterated crossed products is well-known, for example as Proposition~3.11 in \cite{tfb^2}, but we need to know the formula for the induced action with our conventions.

\begin{lemma}\label{cpdecomp}
Suppose that $\alpha:Q\ltimes N\to \Aut A$ is an action of the semidirect product $Q\ltimes N$ by automorphisms of a $C^*$-algebra $A$. Then there is an action $\beta$ of $Q$ on the crossed product $N\ltimes_{\alpha|_N} A$ such that $\beta_p(u_na)=u_{n\cdot p^{-1}}\alpha_{(p,e_N)}(a)$. We write $u$, $v$ and $w$ for the unitary representations of $N$ in $N\ltimes_{\alpha|_N} A$, of $Q$ in $Q\rtimes_\beta (N\ltimes_{\alpha|_N} A)$, and of $Q\ltimes N$ in $(Q\ltimes N)\ltimes_\alpha A$. Then there is an isomorphism of $Q\ltimes_\beta (N\ltimes_{\alpha|_N} A)$ onto $(Q\ltimes N)\ltimes_\alpha A$ which takes $v_p(u_na)$ to $w_{p,n}a$.
\end{lemma}

\begin{proof}
For fixed $p\in Q$, we define $\pi:A\to N\ltimes_{\alpha|_N} A$ by $\pi(a)=i_A(\alpha_{(p,e)}(a))$ and $U:N\to U(N\ltimes_{\alpha|_N} A)$ by $U_n=u_{n\cdot p^{-1}}=i_N(n\cdot p^{-1})$. Then for all $n\in N$, $a\in A$ we have
\begin{align*}
\pi(\alpha_{(e,n)}(a))&=i_A(\alpha_{(p,e)}(\alpha_{(e,n)}(a)))=i(\alpha_{(p,n)}(a))\\
&=i_A(\alpha_{(e,n\cdot p^{-1})}(\alpha_{(p,e)}(a)))=u_{n\cdot p^{-1}}i_A(\alpha_{(p,e)}(a))u^*_{n\cdot p^{-1}}=U_n\pi(a)U_n^*, 
\end{align*}
and the universal property of the crossed product gives a homomorphism $\beta_p:=\pi\ltimes U$ from $N\ltimes_{\alpha|_N} A$ to itself. A straightforward calculation shows that $\beta_p\circ \beta_q=\beta_{pq}$, which implies that each $\beta_p$ is an automorphism, and that $\beta$ is an action of $Q$ on $N\ltimes A$, as claimed. 
 
The map $U:n\mapsto w_{(e,n)}$ is a unitary representation of $N$ in $(Q\ltimes N)\ltimes_\alpha A$, and the pair $(i_A,U)$ is a covariant representation of $(A,N,\alpha|_N)$. Thus there is a homomorphism $i_{N\ltimes A}:=U\ltimes i_A$ of $N\ltimes_{\alpha|_N} A$ into $(Q\ltimes N)\ltimes_\alpha A$. Now $V_p:=w_{(p,e)}$ is a representation of $Q$ in $Q\ltimes_\beta (N\ltimes_{\alpha|_N} A)$ such that $(i_{N\ltimes A},V)$ is covariant, and we get a homomorphism $\pi:=V\ltimes i_{N\ltimes A}$ from $Q\ltimes_\beta (N\ltimes_{\alpha|_N} A)$ to $(Q\ltimes N)\ltimes_\alpha A$. On the other hand, since 
\[
v_pu_{m\cdot p}=\beta_p(u_{m\cdot p})v_p=u_{(m\cdot p)\cdot p^{-1}}v_p=u_mv_p,
\]
there is a unitary representation $W$ of the semidirect product $Q\ltimes N$ in $Q\ltimes_\beta (N\ltimes_{\alpha|_N} xA)$, and this gives a homomorphism back from $(Q\ltimes N)\ltimes_\alpha A$ which is an inverse for $\pi$.
\end{proof}

\section{The multiplicative and Crisp-Laca boundary quotients as crossed products}\label{strucbdaryq}

Recall from Definition~\ref{defn-quotients} that the multiplicative boundary quotient $\partial_{\mult}\TT(\nxxn)$ is the quotient of $\TT(\nxxn)$ obtained by imposing  the relations (T5), and that the Crisp-Laca boundary quotient $\partial\TT(\nxnx)$ is the quotient of $\TT(\nxxn)$ obtained  by imposing  the relations (T5) and (T6). 

In this section we show that $\partial_{\mult}\TT(\nxxn)$ and $\partial\TT(\nxnx)$ are crossed products of an action of the group $\qx$ by automorphisms of the Toeplitz algebra $\TT(\qp)$ and the group algebra $C^*(\Q)$, respectively. 

We start by considering the Toeplitz algebra $\TT(\qp)$. Since the group $(\Q,\qp)$ is totally ordered, every isometric representation of $\qp$ is Nica covariant, and it follows from \cite[Theorem~3.7]{quasilat} (for example) that $\TT(\qp)$ is universal for isometric representations of $\qp$. Indeed we have the following slightly stronger result:
\begin{lemma}\label{faithfulrepsofTT}
Suppose that $W$ is an isometric representation of $\qp$ on a Hilbert space $H$. Then $\pi_W$ is a faithful representation of $\TT(\qp)$ if and only if there exists $r\in \qp$ such that $W_rW_r^*\not=1$.
\end{lemma}
\begin{proof}The statement is a particular case of \cite[Theorem~1]{dou}.
\end{proof}

\begin{proposition}\label{copyofTQ}
There is an isometric representation $W$ of $\qp$ in $\partial_{\mult}\TT(\nxxn)$ such that
\begin{equation}\label{defT}
W_{r}=q_{\mult}(V_aS^bV_a^*)\quad \text{for $r=a^{-1}b$ with $a,b\in \nx$.}
\end{equation}
The induced homomorphism $\pi_W:\TT(\qp)\to \partial_{\mult}\TT(\nxxn)$ is injective.
\end{proposition}

We first check that $W_r$ is well-defined by \eqref{defT}:

\begin{lemma}\label{defTr}
Each $q_{\mult}(V_aS^bV_a^*)$ is an isometry, and 
\begin{equation}\label{welldefrat}
a^{-1}b=c^{-1}d\Longrightarrow q_{\mult}(V_aS^bV_a^*)=q_{\mult}(V_dS^cV_d^*).
\end{equation}
\end{lemma}

\begin{proof}
Since $S$ and $V_a$ are isometries, we have $(V_aS^bV_a^*)^*(V_aS^bV_a^*)=V_aV_a^*$, which implies that $q_{\mult}(V_aS^bV_a^*)$ is an isometry because $q_{\mult}(V_aV_a^*)=1$. For \eqref{welldefrat}, it suffices to assume that $(c,d)=(ak,bk)$ (otherwise swap the two elements). Then we compute using the relation $V_kS^k=SV_k$:
\begin{align*}
V_cS^dV_c^*=V_aV_kS^{bk}V_k^*V_a^*=V_aV_k(S^{k})^bV_k^*V_a^*=V_aS^bV_kV_k^*V_a^*;
\end{align*}
since $q_{\mult}(V_kV^*_k)=1$, \eqref{welldefrat} follows.
 \end{proof}

\begin{proof}[Proof of Proposition~\ref{copyofTQ}]
Lemma~\ref{defTr} implies that that there is a well-defined function $W:\qp\to \partial_{\mult}\TT(\nxxn)$ satisfying \eqref{defT}.
To see that $W$ is a homomorphism, we take $r=a^{-1}b$ and $s=c^{-1}d$, and write $a=a'\gcd(a,c)$, $c=c'\gcd(a,c)$. Then we have
\begin{align}\label{calculationW}
V_aS^bV_a^*V_cS^dV_c^*&=V_aS^bV_{a'}^*V_{c'}S^dV_c^*\\
&=V_aS^bV_{c'}V_{a'}^*S^dV_c^*\quad\quad\ \ \text{by (T3)}\notag\\
&=V_aV_{c'}S^{bc'}S^{a'd}V_{a'}^*V_c^*\quad\text{by (T1)}\notag\\
&=V_{ac'}S^{bc'+a'd}V_{a'c}^*\notag\\
&=V_{ac'}S^{bc'+a'd}V_{ac'}^*,\notag
\end{align}
because $ac'=a'c$. Also
\[
(ac')^{-1}(bc'+a'd)=(ac')^{-1}bc'+(a'c)^{-1}a'd=a^{-1}b+c^{-1}d=r+s.
\]
Applying $q_{\mult}$ to \eqref{calculationW} gives  $W_rW_s=W_{r+s}$.

For each $r=a^{-1}b\in \qp$ we have 
\begin{equation}\label{notone}
q_{\mult}\big(V_a^*(W_rW_r^*)V_a\big)=q_{\mult}\big(V_a^*(V_aS^bS^{*b}V_a^*)V_a\big)=q_{\mult}(S^bS^{*b}).
\end{equation}
The representation $\pi_{\mult}$ of Example~\ref{exbdaryrep} has $\pi_{\mult}(q_{\mult}(S^bS^{*b}))e_x=0$ for $x<b$. Therefore  $q_{\mult}(S^bS^{*b})\not=1$.  Now \eqref{notone} implies that $W_rW_r^*\not=1$. In particular this implies that $W$ is a semigroup of nonunitary isometries, and hence by Lemma~\ref{faithfulrepsofTT} the induced homomorphism $\pi_W$ is injective.
\end{proof}

Since the Toeplitz algebra $(\TT(\Q_+),i_{\qp})=(C^*(T_r:r\in \qp),i_{\qp})$ is universal for isometric representations of $\Q_+$, there is an action $\alpha:\qx\to\Aut \TT(\qp)$ such that 
\begin{equation}\label{alpha-action}
\alpha_a(T_r)=T_{ar} \text{ for $r\in \qp$ and $a\in \qx$}. 
\end{equation}

\begin{prop}\label{idcpToe} Let $W$ be the  isometric representation  of $\qp$ in $\partial_{\mult}\TT(\nxxn)$ of Proposition~\ref{copyofTQ}.
There is a homomorphism $U:\qx\to U\big(\partial_{\mult}\TT(\nxxn)\big)$ such that $U_a=q_{\mult}(V_a^*)$ for $a\in \nx$. Then $(U,\pi_W)$ is a covariant representation of $(\TT(\qp),\qx, \alpha)$, and $U\ltimes\pi_W$ is an isomorphism of $\qx\ltimes_{\alpha}\TT(\qp)$ onto $\partial_{\mult}\TT(\nxxn)$.
\end{prop}

\begin{proof}
The relation (T5) in the multiplicative boundary quotient implies that $a\mapsto q_{\mult}(V_a^*)$ is a homomorphism of the semigroup $\N^\times$ into the unitary group $U\big(\partial_{\mult}\TT(\nxxn)\big)$. Thus it extends uniquely to a unitary representation $U$ of the enveloping group $\qx$.  Let $a,b\in \nx$. The adjoint of relation (T4) implies that
\[
V_a^*S^bV_a=S^{ab}V_a^*V_a.
\]
Hence
\[
U_a\pi_W(T_b)U_a^*=q_{\mult}(V_a^*S^bV_a)=q_{\mult}(S^{ab}V_a^*V_a)=q_{\mult}(S^{ab})=W_{ab}=\pi_W(T_{ab})=\pi_W(\alpha_a(T_b)).
\]
Since the maps $a\mapsto U_a\pi_W(T_b)U_a^*$ and $a\mapsto \pi_W(\alpha_a(T_b))$ agree on $\nx$ and they are both homomorphisms, they also agree on $\qx$. 
Now take $b=c^{-1}d\in \qp$ and $a\in\qx$. Then
\begin{align*}
U_a\pi_W(T_b)U_a^*
&=U_a\pi_W(T_{c^{-1}d})U_a^*=q_{\mult} V_a^*V_cS^dv_c^*a_a)\\
&=U_{ac^{-1}} S^d U_{ac^{-1}}^*=\pi_W(\alpha_{ac^{-1}}(T_d))\\
&=\pi_W(T_{(ac^{-1})d})=\pi_W(\alpha_a(T_{c^{-1}d}))=\pi_W(\alpha_a(T_{b})).
\end{align*}
Thus $(U, \pi_W)$ is a covariant representation.

Since the range of $U\ltimes \pi_W$ is a $C^*$-algebra containing all the generators $q_{\mult}(S)$ and $q_{\mult}(V_a)$, $U\ltimes \pi_W$ is surjective. We know from Proposition~\ref{copyofTQ} that $\pi_W$ is injective. It follows from the presentation (T1)--(T5) that there is  a continuous action $\beta:(\qx)^\wedge\to \Aut\partial_{\mult}\TT(\nxxn)$ such that $\beta_\gamma(q_{\mult}(S))=q_{\mult}(S)$ and 
\[
\beta_\gamma(q_{\mult}(V_a))=\gamma(a)q_{\mult}(V_a)\quad\text{for $\gamma\in (\qx)^\wedge$ and $a\in \N^\times$.}
\]
In other words, we have
\[
(U\ltimes \pi_W)\circ \hat\alpha_\gamma=\beta_\gamma\circ (U\ltimes \pi_W)\quad\text{for $\gamma\in (\qx)^\wedge$.}
\]
Thus it follows from the dual-invariant uniqueness theorem\footnote{This is an old theorem that predates (and motivated) the gauge-invariant uniqueness theorems for Cuntz-Krieger algebras, graph algebras, Cuntz-Pimsner algebras, etc. A sledgehammer reference is \cite[Corollary~4.3]{R}, which is about crossed products by coactions but specialises to what we want for abelian groups. But the case of a discrete abelian group is much easier.} that $U\ltimes \pi_W$ is faithful. 
\end{proof}

We can use the isomorphism of Proposition~\ref{idcpToe} to study the structure of the multiplicative boundary quotient by analysing the  crossed product $\qx\ltimes_{\alpha}\TT(\qp)$, and we will do this in the next section. However, we can already get quite a bit more information using Douglas' work on the Toeplitz algebras of dense subgroups of $\R$ \cite{dou}. The canonical representation $u$ of $\qp$ in the group algebra $C^*(\Q)$ induces a homomorphism $\pi_u:\TT(\qp)\to C^*(\Q)$. Douglas proved that the kernel of $\pi_u$ is the commutator ideal $\CC(\qp)$ of $\TT(\qp)$, and that the ideal $\CC(\qp)$ is simple (see the Corollary on \cite[page~147]{dou}). Thus there is an exact sequence 
\[
0\longrightarrow \CC(\qp)\longrightarrow \TT(\qp)\stackrel{\pi_u}{\longrightarrow}C^*(\Q)\longrightarrow 0.
\]
For $a\in\qx$, the automorphisms $\alpha_a$ map commutators to commutators, and so the ideal $\CC(\qp)$ is invariant for the action $\alpha$ defined at \eqref{alpha-action}. Thus it follows from \cite[Proposition~3.19]{tfb^2}, for example, that there is an exact  sequence 
\begin{equation}\label{cpexact}
0\longrightarrow \qx\ltimes_\alpha\CC(\qp)\longrightarrow \qx\ltimes_\alpha\TT(\qp)\stackrel{\id\ltimes \pi_u}{\longrightarrow} \qx\ltimes_\alpha C^*(\Q)\longrightarrow 0.
\end{equation}
Since the quotient $\qx\ltimes_\alpha C^*(\Q)\cong\qx\ltimes_\alpha C(\widehat \Q)$ at \eqref{cpexact} is the $C^*$-algebra of a transformation-group, we can compute its ideal structure by studying the action of $\qx$ on $\widehat{\Q}$, and we do this in the next section.

We finish this section by identifying the image of   the ideal $\qx\ltimes_\alpha\CC(\qp)$ in the multiplicative boundary quotient and the image of the quotient  by $\qx\ltimes_\alpha\CC(\qp)$ as the boundary quotient.

\begin{prop}\label{prop-commutingdiagram-revised}
\begin{enumerate}
\item\label{prop-commutingdiagram1} 
There is a homomorphism $\widetilde U:\qx\to U\big(\partial\TT(\nxxn)\big)$ such that $\widetilde U_a=q_{\bdry}(V_a^*)$ for $a\in\nx$.
\item\label{prop-commutingdiagram2}  There is a unitary representation $\widetilde W: \Q\to \partial\TT(\nxxn)$ such that $\widetilde W_{a^{-1}b}=q_{\bdry}(V_aS^bV_a^*)$ for $a,b\in\nx$.  
\item\label{prop-commutingdiagram3}  $(\widetilde U,\pi_{\widetilde W})$ is a covariant representation of $(C^*(\Q),\qx,\alpha)$. 
\item\label{prop-commutingdiagram4} 
Let $q:\partial_{\mult}\TT(\nxxn)\to \partial\TT(\nxxn)$ be the quotient map. Then $\widetilde U\ltimes\pi_{\widetilde W}$ is an isomorphism of $\qx\ltimes_{\alpha}C^*(\Q)$ onto $\partial\TT(\nxxn)$ such that
\begin{equation}\label{commdiagidealbdryquotient}
\xymatrix{
0\longrightarrow \qx\ltimes_\alpha\CC(\qp) \ar[r]\ar[d]_{(U\ltimes\pi_W)|}
&\qx\ltimes_\alpha \TT(\Q_+)\ar[r]^{\id\ltimes \pi_u\quad}\ar[d]_{U\ltimes\pi_W}
&\qx\ltimes_\alpha C^*(\Q)\ar[d]^{\widetilde U\ltimes \pi_{\widetilde W}}\longrightarrow 0
\\
0\longrightarrow \ker q\ar[r]
&\partial_{\mult}\TT(\nxxn)\ar[r]^{q}
&\partial\TT(\nxxn)\longrightarrow 0}
\end{equation}
commutes. In particular, $(U\ltimes\pi_W)(\qx\ltimes_\alpha\CC(\qp))=\ker q$.
\end{enumerate}
\end{prop}

\begin{proof} By Proposition~\ref{idcpToe}, there is a homomorphism $U:\qx\to U\big(\partial_{\mult}\TT(\nxxn)\big)$ such that $U_a=q_{\mult}(V_a^*)$ for $a\in \nx$.  Since $q_{\bdry}=q\circ q_{\mult}$ it follows that  $\widetilde U=q\circ U$ has the required properties.   This gives (\ref{prop-commutingdiagram1}). 

By Proposition~\ref{copyofTQ}, there is an isometric representation $W$ of $\qp$ in $\partial_{\mult}\TT(\nxxn)$ such that
$
W_{r}=q_{\mult}(V_aS^bV_a^*)\quad \text{for $r=a^{-1}b$ with $a,b\in \nx$.}
$
For $r\in \Q$ we  set 
\[
\widetilde W_r=\begin{cases}
q(W_r)&\text{if $r>0$}\\
1&\text{if $r=0$}\\
q(W_{-r}^*)&\text{if $r<0$}.
\end{cases}
\]
It is straightforward to check that (T5) and (T6) imply that each $\widetilde W_r$ is unitary: for example, if $r=a^{-1}b$ with $a,b\in\nx$, then
\[
\widetilde W_r\widetilde W_{-r}=q_{\bdry}\big(V_aS^bV_a^*V_a(S^{-b})^*V_a^*)\big)=1
\]
because $q_{\bdry}(SS^*)=1$ and $q_{\bdry}(V_aV_a^*)=1$.  

To check that $\widetilde W$ is a representation, fix $r, s\in \Q$. If either both $r,s\geq 0$ or both $r,s\leq 0$, then it is straightforward to check that $\widetilde W_{r+s}=\widetilde W_r\widetilde W_s$. So take $r=a^{-1}b>0$, $s=c^{-1}d<0$ where $a, b, c\in\nx$, and $d\in-\nx$. Write $a=a'\gcd(a,c)$ and $c=c'\gcd(a,c)$ where $a',c'\in\nx$, and calculate as at \eqref{calculationW} to get that
\begin{align*}
\widetilde W_r\widetilde W_s&=q_{\bdry}(V_aS^bV_a^*)q_{\bdry}(V_cS^{-d}V_c^*)^*=q_{\bdry}(V_aS^bV_{a'}^*V_{c'}(S^{-d})^*V_c^*)
\\
&=q_{\bdry}(V_aS^bV_{c'}V_{a'}^*(S^{-d})^*V_c^*)=q_{\bdry}(V_aV_{c'}S^{bc'}(S^{-a'd})^*V_{a'}^*V_c^*)
\\
&=q_{\bdry}(V_{ac'}S^{bc'}(S^{-a'd})^*V_{a'c}^*).
\end{align*}
Now we note that $r+s=\frac{bc'+a'd}{a'c}$ and distinguish two cases. First, suppose that  $r+s\geq 0$. Then $bc'+a'd\geq 0$ and using $q_{\bdry}(SS^*)=1$ we get that 
\[
\widetilde W_{r+s}=q_{\bdry}(V_{a'c}S^{bc'+ad'}V_{a'c}^*)=q_{\bdry}(V_{a'c}S^{bc'}(S^{-ad'})^*V_{a'c}^*)=\widetilde W_r\widetilde W_s.
\] Second,  suppose that  $r+s<0$. Then $bc'+a'd< 0$ and again
\[\widetilde W_{r+s}=q_{\bdry}(V_{a'c}(S^{-(bc'+ad')})^*V_{a'c}^*)=q_{\bdry}(V_{a'c}S^{bc'}(S^{-ad'})^*V_{a'c}^*)=\widetilde W_r\widetilde W_s.\] So $\widetilde W$ is a unitary representation of $\Q$ in $\partial\TT(\nxxn)$, which is   (\ref{prop-commutingdiagram2}). 

To see $(\widetilde U, \pi_{\widetilde W})$ is covariant for $(\qx, C^*(\Q), \alpha)$, let $a\in \qx$ and $r\in \Q$. We have $\pi_{\widetilde W}(\alpha_a(u_r))=\pi_{\widetilde W}(u_{ar})=\widetilde W_{ar}$. Using the covariance of $(U,\pi_W)$ we have
\begin{align*}
\widetilde U_a\pi_{\widetilde W}(u_s)\widetilde U_a^*
&=\widetilde U_a \widetilde W_s\widetilde U_a^*\\
&=\begin{cases}
q(U_aW_rU_a^*) &\text{if $r>0$}\\
q(U_aU_a^*) &\text{if $r=0$}\\
q(U_aW_{-r}^*U_a^*) &\text{if $r<0$}\\
\end{cases}\\
&=\begin{cases}
q(W_{ar}) &\text{if $r>0$}\\
1 &\text{if $r=0$}\\
q(W_{-ar}^*) &\text{if $r<0$}\\
\end{cases}\\
&=\widetilde W_{ar}.
\end{align*}
 This gives (\ref{prop-commutingdiagram3}).

For (\ref{prop-commutingdiagram4}) we start by proving that the right-hand square of the diagram \eqref{commdiagidealbdryquotient} commutes. It suffices to check this for the generators of $\TT(\qp)$ and $\qx$. Let $r\in\qp$. Then the upper route is
$
T_r\mapsto u_r\mapsto \widetilde W_r
$
and the lower route is $T_r\mapsto W_r\mapsto \widetilde W_r$. Next let $a\in\qx$. Then the upper route is $a\mapsto a \mapsto \widetilde U_a$ and the lower route is 
$a\mapsto U_a \mapsto \widetilde U_a$. It follows that the right-hand square commutes.

Next we  show that  $\widetilde U\ltimes \pi_{\widetilde W}$ is an isomorphism. To do this, we view $\qx\ltimes_\alpha C^*(\Q)$ as the universal $C^*$-algebra generated by unitary representations $X$ of $\qx$  and $Y$ of $\Q$ in $\qx\ltimes_\alpha C^*(\Q))$ satisfying $X_aY_rX_a^*=Y_{ar}$.  Then it is straightforward to verify that  $\{Y_1, X_{a^{-1}}:a\in \nx\}$ satisfies the relations (T0--T6), or equivalently, that $L(a,m)=X_{a^{-1}}Y_m$ is a Nica-covariant representation of $\nxxn$ such that $\pi_L:\TT(\nxxn)\to \qx\ltimes_\alpha C^*(\Q)$ factors through the quotient $\partial\TT(\nxxn)$. Let $\widetilde\pi_L:\partial\TT(\nxxn)\to \qx\ltimes_\alpha C^*(\Q)$ be the representation such that $\pi_L=\widetilde\pi_L\circ q_{\bdry}$.  We claim that $\widetilde\pi_L$ is an inverse for $\widetilde U\ltimes \pi_{\widetilde W}$. We check on generators. For $(a,m)\in\nxxn$ we have
\[
(\widetilde U\ltimes \pi_{\widetilde W})\circ\widetilde\pi_L(q_{\bdry}(T_{(a,m)})
=\widetilde U\ltimes \pi_{\widetilde W}(X_{a^{-1}}Y_m)
=\widetilde U_{a^{-1}}\widetilde W_m
=q_{\bdry}(V_a)q_{\bdry}(S^m)
=q_{\bdry}(T_{(a,m)}),
\]
and it follows that $(\widetilde U\ltimes \pi_{\widetilde W})\circ\widetilde\pi_L=\id$.  For $a,b\in \nx$ we have
\[
\widetilde\pi_L\circ(\widetilde U\ltimes \pi_{\widetilde W})(X_{a^{-1}b})=\widetilde\pi_L(\widetilde U(X_{a^{-1}b}))=\widetilde\pi_L(\widetilde U_{a^{-1}b})
=\widetilde\pi_L (q_{\bdry}V_aV_b^*)=X_{a^{-1}}X_{b^{-1}}^*=X_{a^{-1}b}.
\]
For $r\in\Q$ we have
\begin{align*}
\widetilde\pi_L\circ(\widetilde U\ltimes \pi_{\widetilde W})(Y_r)=\widetilde\pi_L(\widetilde W_r)
&=\begin{cases}
\widetilde\pi_L q_{\bdry}(W_{a^{-1}b})&\text{if $r=a^{-1}b>0$}\\
1&\text{if $r=0$}\\
\widetilde\pi_L q_{\bdry}(W_{-a^{-1}b})^*&\text{if $r=a^{-1}b<0$}
\end{cases}\\
&=\begin{cases}
L (a,b)L(a, 0)^*&\text{if $r=a^{-1}b>0$}\\
1&\text{if $r=0$}\\
(L (a,-b)L(a, 0)^*)^*&\text{if $r=a^{-1}b<0$}
\end{cases}\\
&=\begin{cases}
X_{a^{-1}}Y_bX_{a^{-1}}^*&\text{if $r=a^{-1}b>0$}\\
1&\text{if $r=0$}\\
(X_{a^{-1}}Y_{-b}X_{a^{-1}}^*)^*&\text{if $r=a^{-1}b<0$}
\end{cases}\\
&=\begin{cases}
 Y_{a^{-1}b}&\text{if $r=a^{-1}b>0$}\\
Y_0&\text{if $r=0$}\\
Y_{-a^{-1}b}^*&\text{if $r=a^{-1}b<0$}
\end{cases}
\end{align*}
which is just $Y_r$. Thus $\widetilde\pi_L\circ(\widetilde U\ltimes \pi_{\widetilde W})=\id$, and $\widetilde U\ltimes \pi_{\widetilde W}$ is an isomorphism.

Next we show that $\widetilde U\ltimes\pi_{\widetilde W}\big(\qx\ltimes_\alpha\CC(\qp)\big)=\ker q_{\bdry}$.   Let $x\in  \qx\ltimes_\alpha\CC(\qp)$. Then $0=\pi_u(x)$ and hence $0=(\widetilde U\ltimes\widetilde W)\circ \pi_u(x)=q_{\bdry}\circ (U\ltimes W)(x)$. Thus $(U\ltimes W)(x)\in\ker q_{\bdry}$. Thus $\widetilde U\ltimes\pi_{\widetilde W}\big(\qx\ltimes_\alpha\CC(\qp)\big)\subset\ker q_{\bdry}$, and the other inclusion follows from a similar argument using the inverses of $\widetilde U\ltimes\pi_{\widetilde W}$ and $U\ltimes W$.  
Thus $\widetilde U\ltimes\pi_{\widetilde W}\big(\qx\ltimes_\alpha\CC(\qp)\big)=\ker q_{\bdry}$, and it is now clear the left-hand square commutes as well.
\end{proof}

\section{The structure of the crossed product $\qx\ltimes_{\alpha}\TT(\qp)$.}

In this section we study the structure of the crossed products in the short exact sequence
\begin{equation*}
0\longrightarrow \qx\ltimes_\alpha\CC(\qp)\longrightarrow \qx\ltimes_\alpha\TT(\qp)\stackrel{\id\ltimes \pi_u}{\longrightarrow} \qx\ltimes_\alpha C^*(\Q)\longrightarrow 0.
\end{equation*}
 The main theorem of this section, Theorem~\ref{primcp},  describes the primitive-ideal space of  $\qx\ltimes_\alpha\TT(\qp)$.  At the end of the section we pull this result across to $\partial_{\mult}\TT(\nxxn)$ using the isomorphism of Proposition~\ref{idcpToe}.

\begin{theorem}\label{primcp}
Consider the action $\alpha:\qx\to\Aut\TT(\qp)$ such that $\alpha_a(T_r)=T_{ar}$. Let $\psi$ be the homomorphism of $\TT(\qp)$ onto $\C$ such that $\psi(T_r)=1$ for $r\in \qp$. Then there is a homomorphism $\id\ltimes\psi$ of $\qx\ltimes_{\alpha}\TT(\qp)$ onto $C^*(\qx)=\qx\ltimes_{\id}\C$. For $\gamma\in (\qx)^\wedge$,  let $\epsilon_\gamma$ be the character of $C^*(\qx)$ such that $\epsilon_\gamma(u_a)=\gamma(a)$. 
\begin{enumerate}
\item The map $\gamma\mapsto \ker\big(\epsilon_\gamma\circ (\id\ltimes\psi)\big)$ is a homeomorphism of the compact space $(\qx)^\wedge$ onto a closed subset $F$ of $\Prim\big(\qx\ltimes_{\alpha}\TT(\qp)\big)$; the ideal of $\qx\ltimes_{\alpha}\TT(\qp)$ with hull $F$ is $\ker(\id\ltimes\psi)$. 
\item The other primitive ideals are $O:=\{0\}$ and $I:=\qx\ltimes_\alpha\CC(\qp)$. \item The set $\{O\}$ is open and dense in $\Prim\big(\qx\ltimes_{\alpha}\TT(\qp)\big)$, and the set $\{I\}$ is open and dense in $\Prim\big(\qx\ltimes_{\alpha}\TT(\qp)\big)\setminus \{O\}$.
\item Both $I$ and $\ker(\id\ltimes\psi)/I$ are simple.
\end{enumerate}
\end{theorem}

To study $\TT(\Q_+)$ we  realise it as a corner in a crossed product, which goes back at least as far as Murphy \cite{M}, and was used extensively by Phillips and Raeburn in \cite{PR}, which will be our main reference. We define
\[
B_{\Q}(-\infty,\infty]:=\clsp\big\{1_{[r,\infty)}:r\in \Q\big\}\subset \ell^\infty(\R).
\]
Then $B_{\Q}(-\infty,\infty]$ is a $C^*$-subalgebra of $\ell^\infty(\R)$, and consists of the functions $f:\R\to \C$ that are continuous at all $x\in \R\setminus\Q$, are right-continuous at $r\in\Q$,  have limits as $x\to r^-$ and $x\to r^+$ for $r\in \Q$ (not necessarily equal), and have limits as $x\to \infty$ (see \cite[Proposition~3.1]{PR}). We consider also the $C^*$-
subalgebra
\[
B_{\Q}(\R)=\big\{f\in B_{\Q}(-\infty,\infty]:f(x)\to 0\text{ as $x\to \infty$}\big\}.
\]
The map $\epsilon_{\infty}:f\mapsto \lim_{x\to\infty}f(x)$ is a homomorphism, and there is  an exact sequence
\[
0\longrightarrow B_{\Q}(\R)\longrightarrow B_{\Q}(-\infty,\infty]\stackrel{\epsilon_\infty}{\longrightarrow}\C\longrightarrow 0.
\]
We define $\tau:\Q\to \Aut B_{\Q}(-\infty,\infty]$ by $\tau_r(f)(x)=f(x-r)$. The ideal $\bqr$ is invariant for the action $\tau$, and hence \cite[Proposition~3.19]{tfb^2} gives an exact sequence
\[
0\longrightarrow \Q\ltimes_\tau B_{\Q}(\R)\longrightarrow \Q\ltimes_\tau B_{\Q}(-\infty,\infty]\stackrel{\id\ltimes\epsilon_\infty}{\longrightarrow}\Q\ltimes\C
=C^*(\Q)\longrightarrow 0.
\]
With $P:=i_{B_\Q}(1_{[0,\infty)})$, there is an isomorphism $\pi$ of $\TT(\Q_+)$ onto the corner $P\big(\Q\ltimes_\tau B_{\Q}(-\infty,\infty]\big)P$ such that $\pi(T_r)=\pi_{\Q}(r)P$ for $r\in \qp$. (This follows, for example, by applying \cite[Proposition~3.2]{PR} to the trivial cocycle $\sigma\equiv 1$.) This isomorphism carries the commutator ideal $\CC(\qp)$ onto the corner $P\big(\Q\ltimes_{\tau} B_{\Q}(\R)\big)P$. 

We now want to study the action of $\qx$ on the corner that corresponds to the action $\alpha$ on the Toeplitz algebra. Our key observation is that there is a related action of the semidirect product $\qxxq$ on $B_{\Q}(-\infty,\infty]$. Indeed, we define $\gamma:\qxxq\to \Aut B_{\Q}(-\infty,\infty]$ by $\gamma_{a,r}(f)(x)=f(ax-r)$. Observe that we have $\gamma_{0,r}=\tau_r$.  Then Lemma~\ref{cpdecomp} gives:

\begin{lemma}
There is an action 
\[
\beta:\qx\to \Aut\big(\qp\ltimes_\tau B_{\Q}(-\infty,\infty]\big)\text{ such that }\beta_a(u_rf)=u_{ra^{-1}}\gamma_{a,0}(f).
\] 
We write $u$, $v$ and $w$ for the canonical unitary representations of $\qp$, $\qx$ and $\qxxq$ in $\qxxq$ in the respective crossed crossed products by $\tau$, $\beta$ and $\gamma$.  Then there is an
isomorphism\footnote{ It is possibly important that the isomorphism $U\ltimes \pi_W$ of Proposition~\ref{idcpToe} does not carry the action $\alpha$ into $\beta$. The unitary representation  $U$ in that Proposition carries $a$ into $q_{\mult}(V_a^*)$. The action $\beta$ comes from the action $\gamma$ of $\qx$ rather than $\alpha$.}
\[
\pi:\qx\ltimes_\beta\big(\Q\ltimes_\tau B_{\Q}(-\infty,\infty]\big)\longrightarrow(\qxxq)\ltimes_\gamma B_{\Q}(-\infty,\infty]
\] 
such that (in the notation of \S\ref{cps}) $\pi(v_au_rf)=w_{a,r}f$.
\end{lemma}

The action $\beta$ of $\qx$ fixes the projection $P$, and hence restricts to an action  on $P\big(\Q\ltimes_\tau B_{\Q}(-\infty,\infty]\big)P$. When we view this corner as the Toeplitz algebra $\TT(\qp)$, the generating isometries are $T_r=Pu_rP$. Thus $\beta$ is characterised by the formula $\beta_a(T_r)=T_{a^{-1}r}$, and is not quite the action $\alpha$ which we have been studying. But it is very closely related, and the crossed product is the same. More precisely, we have:

\begin{lemma}\label{addinv}
Suppose that $(A,G,\alpha)$ is a dynamical system in which the group $G$ is abelian. Define $\beta:G\to\Aut A$ by $\beta_g=\alpha_g^{-1}$. Then $(G\ltimes_{\alpha} A,i_A, i_G\circ{(g\mapsto g^{-1})})$ is a crossed product for $(A,G,\beta)$, and  $G\ltimes_{\beta} A$ is canonically isomorphic to $G\ltimes_{\alpha}A$.
\end{lemma}

\begin{proof}
Since $G$ is abelian, the map $\inver:g\mapsto g^{-1}$ is an isomorphism of $G$ onto $G$. Now a calculation shows that if $\pi:A\to B(H)$ and $U:G\to U(H)$, then $(\pi,U)$ is a covariant representation of $(A,G,\alpha)$ if and only if $(\pi,U\circ{\inver})$ is a covariant representation of $(A,G,\beta)$. So $(G\ltimes_{\alpha} A,i_A, i_G\circ{\inver})$ is a crossed product for $(A,G,\beta)$. 
\end{proof}

The isomorphism of $\TT(\qp)$  onto the corner in $\Q\ltimes_\tau B_{\Q}(-\infty,\infty]$ restricts to an isomorphism of the commutator ideal $\CC(\qp)$ onto a corner in $\Q\ltimes_\tau B_{\Q}(\R)$, and we can realise $\qx\ltimes_{\beta}\CC(\qp)$ as a corner in $\qx\ltimes_\beta\big(\Q\ltimes_\tau B_{\Q}(\R)\big)$.

The automorphisms $\big\{\gamma_{a,r}:(a,r)\in \qxxq\big\}$ satisfy $\epsilon_\infty(\gamma_{a,r}(f))=\epsilon_\infty(f)$, and hence induce automorphisms of the ideal $B_{\Q}(\R)=\ker\epsilon_\infty$. The following technical result may be of some independent interest --- certainly, we were surprised to discover it.

\begin{prop}\label{cpsimple2}
The $C^*$-algebra $(\qxxq)\ltimes_{\gamma}B_{\Q}(\R)$ is simple.
\end{prop}

We will prove this using a result of  Archbold and Spielberg \cite{AS}. To apply their result, we have to prove that the action of $\qxxq$ on the spectrum of $B_{\Q}(\R)$ is topologically free and minimal. We will use the following criterion for topological freeness.

\begin{lemma}\label{crittopfree}
Suppose that $G$ acts on a locally compact space $X$, and that $\big\{x\in X:S_x=\{e\}\big\}$ is dense in $X$. Then the action on $X$ is topologically free in the sense of \cite[Definition~1]{AS}.
\end{lemma}

\begin{proof}
Suppose that $g_1,\cdots,g_n\in G\setminus\{e\}$. We want to show that that
\begin{equation}\label{finiteinter}
\bigcap_{i=1}^n\big\{x\in X:g_i\cdot x\not= x\big\}
\end{equation}
is dense in $X$. Suppose that $U$ is an open subset of $X$. Then there exists $y\in U$ such that $S_y=\{e\}$. Then $g_i\cdot y\not=y$ for $1\leq i\leq n$, and hence $y\in U$ belongs to the intersection \eqref{finiteinter}.
\end{proof}

The spectrum of $\bqr$ is described in Lemma~3.6 of \cite{PR}. As a set
\begin{equation*}
\bqr^\wedge=\{\epsilon_x:x\in\R\}\cup\{\epsilon_\lambda^-:\lambda\in\Q\}
\end{equation*} 
where $\epsilon_x$ is evaluation at $x$ and $\epsilon_\lambda^-(f)=\lim_{x\to\lambda^-}f(x)$ for $f\in \bqr$. Lemma~3.6 of \cite{PR} also describes the topology of $\bqr^\wedge$  in terms of convergence of sequences.

\begin{lemma}\label{lem-bqr^action} Let $\gamma:\qx\ltimes\qp\to\Aut\bqr$ be the action defined by $\gamma_{a,r}(f)(x)=f(ax-m)$. For $x\in\R$ and $\lambda\in\Q$, the induced action of $\qx\ltimes\qp$ on $\bqr^\wedge$ is given by
\begin{equation}\label{induced_action}
(a,r)\cdot \epsilon_x=\epsilon_{a^{-1}x+a^{-1}r}\quad\text{and}\quad(a,r)\cdot \epsilon_\lambda^-=\epsilon^-_{a^{-1}\lambda+a^{-1}r}.
\end{equation} 
This action is topologically free and minimal.
\end{lemma}

\begin{proof}
Let $x\in\R$, $\lambda\in\Q$, $(a,r)\in\qx\ltimes\qp$ and $f\in\bqr$. Then 
\begin{align*}
((a,r)\cdot \epsilon_x)(f)&
=\epsilon_x\circ \gamma_{a,r}^{-1}(f)
=\gamma_{a^{-1},-a^{-1}r)}(f)(x)\\
&=f(a^{-1}x+a^{-1}r)=\epsilon_{a^{-1}x+a^{-1}r}(f), 
\end{align*}
and, similarly,
\begin{align*}
((a,r)\cdot \epsilon_\lambda^-)(f)
&=\lim_{x\to\lambda^-} \gamma_{a,r}^{-1}(f)(x)
=\lim_{x\to\lambda^-} f(a^{-1}x+a^{-1}r)\\
&=\lim_{x\to(a^{-1}\lambda+a^{-1}r)^-} f(x)
=\epsilon^-_{a^{-1}\lambda+a^{-1}r}(f).
\end{align*}
This proves that the induced action is given by \eqref{induced_action}.

We now verify that the action is topologically free. Let $x\in\R$. Then 
\[
(a,r)\cdot\epsilon_x=\epsilon_x\Longleftrightarrow r=(a-1)x.\]
Thus if $x\in\R\setminus\Q$, then $S_{\epsilon_x}=\{(1,0)\}$.  We will show that $\{\epsilon_x:x\in\R\setminus\Q\}$ is dense in $\bqr^\wedge$. Let $\lambda\in \Q$, and consider $\epsilon_\lambda$ and $\epsilon_\lambda^-$. Choose $\{x_i\}, \{y_i\}\subset \R\setminus \Q$ such that $x_i\to\lambda^+$ and $y_i\to\lambda^-$. Then  $\epsilon_{x_i}\to\epsilon_\lambda$ and $\epsilon_{y_i}\to \epsilon_\lambda^-$ by Lemma~3.6 of \cite{PR}. Thus $\{\epsilon_x:x\in\R\setminus\Q\}$ is dense in $\bqr^\wedge$. Now Lemma~\ref{crittopfree} implies that the action is topologically free. 

The action is minimal when  $\bqr$ has no nontrivial  $(\qx\ltimes\qp)$-invariant ideals. Thus, since $\bqr$ is commutative, the action is minimal if and only if the only non-empty closed $(\qx\ltimes\qp)$-invariant subset is $\bqr^\wedge$. Thus it suffices to show that every orbit closure is $\bqr^\wedge$. Let $\phi\in \bqr^\wedge$. First  suppose that $\phi=\epsilon_x$ for some $x\in\R$. Let $y\in\R$.  Choose sequences $\{s_i\}, \{r_i\}\subset\Q$ such that $x+s_i\to y^+$ and $x+r_i\to y^-$. Then by Lemma~3.6 of \cite{PR},
\[
(1, s_i)\cdot\phi=\epsilon_{x+s_i}\to\epsilon_y\quad\text{and}\quad (1, r_i)\cdot\phi=\epsilon_{x+r_i}\to
\begin{cases}
\epsilon_y &\text{if $y\notin \Q$}\\
\epsilon_y^- &\text{if $y\in \Q$.}
\end{cases}
\]
It follows that $\overline{(\qx\ltimes\qp)\cdot\phi}=\bqr^\wedge$. Second, suppose that $\phi=\epsilon_\lambda^-$ for some $\lambda\in\Q$. Then $(\qx\ltimes\qp)\cdot\phi=\{\epsilon_r^-:r\in\Q\}$. Let $y\in\R$. Choose $\{r_i\}\subset \Q$ such that $r_i\to y^+$. For each $r_i$ and $j\geq 1$, choose $z_{i,j}>y$ such that $r_i-\frac{1}{j}<z_{i,j}<r_i$. Then $z_{i,j}\to r_i^-$ as $j\to\infty$. For $f\in \bqr$ the right continuity of $f$ gives\[
\lim_{i\to\infty} \epsilon_{r_i}^-(f)=\lim_{i\to\infty}\lim_{x\to r_i^-}f(x)=\lim_{i\to\infty}\lim_{j\to\infty}f(z_{i,j})=f(y)=\epsilon_y(f).
\] 
using the right continuity of $f$.
Thus again we have $\overline{(\qx\ltimes\qp)\cdot\phi}=\bqr^\wedge$, and the action is minimal. 
\end{proof}

\begin{proof}[Proof of Proposition~\ref{cpsimple2}]
Lemma~\ref{lem-bqr^action} says that the action of $\qxxq$ on $B_{\Q}(\R)^\wedge$ satisfies the first two hypotheses of the Corollary to Theorem~2 in \cite{AS}. The other hypothesis (``regularity'') is automatic because the group  $\qxxq$ is amenable (see the comment before the Corollary). Thus the result follows from that Corollary.\end{proof}

Since $\qx\ltimes_{\alpha}\CC(\qp)$ is isomorphic to $\qx\ltimes_{\beta}\CC(\qp)$, and since the latter is isomorphic to a corner in the simple algebra 
\[\qx\ltimes_{\beta}\big(\Q\ltimes_\tau B_{\Q}(\R)\big)\cong (\qx\ltimes\Q)\ltimes_\gamma B_{\Q}(\R) \] we can now use Proposition~\ref{cpsimple2} to prove following corollary.

\begin{cor}\label{cpsimple}
The $C^*$-algebra $\qx\ltimes_{\alpha}\CC(\qp)$ is simple.
\end{cor}

\begin{proof}
The Toeplitz algebra $\TT(\qp)$ is the subalgebra of $\qp\ltimes_\tau B_{\Q}(-\infty,\infty]$ generated by the operators $Pu_mP= 1_{[0,\infty]}u_m1_{[0,\infty]}$. Since 
\[
\gamma_{(a,0)}(1_{[0,\infty]})(x)=1_{[0,\infty]}(ax)=\begin{cases}1&\text{if $ax\geq 0\Longleftrightarrow x\geq 0$}\\
0&\text{if $ax<0\Longleftrightarrow x<0$,}
\end{cases} 
\]
we have $\gamma_{(a,0)}(1_{[0,\infty]})=1_{[0,\infty]}$. Thus $\beta_a$ is the automorphism of the Toeplitz algebra $\TT(\qp)=C^*(T_p)$ such that $\beta_a(T_p)=T_{a^{-1}p}$. We know from Corollary~\ref{cpsimple2} that $(\qxxq)\ltimes_\gamma B_{\Q}(\R)$ is simple, and hence so is $\qx\ltimes_\beta\big(\qp\ltimes_\tau B_{\Q}(\R)\big)$. Thus so is the corner associated to the projection $1_{[0,\infty]}$, which is $\qx\ltimes_\beta \CC(\qp)$.

Since the group $\qx$ is abelian, the crossed product by the action $\alpha:a\mapsto \beta_a^{-1}$ is isomorphic to $\qx\ltimes_{\beta}\CC(\qp)$ (see Lemma~\ref{addinv}). Thus $\qx\ltimes_{\alpha}\CC(\qp)$ is simple. 
\end{proof}

We now analyse the structure of the crossed product $\qx\ltimes_\alpha C^*(\Q)$ by the action induced by the action $\alpha:\qx\to \Aut\TT(\qp)$. We use the Fourier transform $\F$ to identify $C^*(\Q)$ with $C(\widehat\Q)$: by convention, we use the transform such that $\F(u_r)(\gamma)=\gamma(r)$, and then the action $\alpha$ is given by $\alpha_a(\gamma)(r)=\gamma(ar)=(a\cdot\gamma)(r)$.

Next we realise  $\widehat \Q$ as a solenoid. Let $m,n\in \nx$. Then  $m\mid n$ implies that  $n=mc$ for some $c\in \nx$, and $m^{-1}\Z\subset n^{-1}\Z$ with inclusion map $m^{-1}k\mapsto n^{-1}ck$. We then view $\Q$ as the direct limit
\[
\Q=\bigcup_{n\in\nx}n^{-1}\Z=\varinjlim_{n\in\nx}n^{-1}\Z.
\]
The dual group of each $n^{-1}\Z$ is $\T$, and hence we can view $\widehat \Q$ as the inverse limit
\[
\mathcal{S}=\varprojlim _{n\in\nx}\T=\Big\{(z_n): z_n\in\T \text{\ and\ }n=mc\Rightarrow z_m=z_n^c\Big\}.
\]
For $z=(z_n)\in\mathcal{S}$ we have the character $\gamma_z$ in $\widehat \Q$ with formula $\gamma_z(n^{-1}k)=z_n^k$. Let $a,b\in\nx$.
The action of $\qx$ is given by $ab^{-1}\cdot\gamma_z(n^{-1}k)=\gamma_n(ab^{-1}n^{-1}k)=\gamma_z((nb)^{-1}ak)=z_{nb}^{ak}=\gamma_{z^a}((nb)^{-1}k)$.  Thus $(ab^{-1}\cdot z)$ is the sequence $z^{ab^{-1}}$ in $\mathcal{S}$ with
\[
(z^{ab^{-1}})_n=z_{nb}^a.
\]
 
The constant sequence $1$ is fixed under the action of $\qx$. Hence the evaluation map $\epsilon_1:f\mapsto f(1)$ gives an exact sequence
\begin{equation*}
0\longrightarrow C_0(\mathcal{S}\setminus \{1\})\longrightarrow C(\mathcal{S})\stackrel{\epsilon_1}{\longrightarrow}\C\longrightarrow 0
\end{equation*}
which is preserved by the action of $\qx$. (We observe that the homomorphism $\epsilon_1$ is essentially the homomorphism $\phi:C^*(\Q)\to\C$ such that $\phi(u_r)=1$ for all $r$, in the sense that $\phi=\epsilon_1\circ\F$. This homomorphism is used  in Proposition~\ref{primofC*(qxxq)} below, and satisfies $\phi\circ\pi_u=\psi$.) Thus it follows from \cite[Proposition~3.19]{tfb^2} that there is an exact sequence
\begin{equation}\label{cpseqexact}
0\longrightarrow \qx\ltimes C_0(\mathcal{S}\setminus \{1\})\longrightarrow \qx\ltimes C(\mathcal{S})\stackrel{\id\ltimes\epsilon_1}{\longrightarrow}\qx\ltimes \C=C^*(\qx)\longrightarrow 0.
\end{equation}

\begin{lemma}\label{free+min}
The action of $\qx$ on $\mathcal{S}\setminus\{1\}$ is free and minimal. 
\end{lemma}
 
\begin{proof}
For freeness, we take $z=(z_n)\in\mathcal{S}$. We first suppose that $a\cdot z=z$ and aim to show that $z=1$. We have $z_n^{a}=z_n$ for all $n$, and hence $z_n^{a-1}=1$ for all $n$. Then $1=z_{n(a-1)}^{a-1}=z_n$ for all $n$, and $z=1$, as we wanted.  

Next we take $a,b\in\nx$ such that $a> b$ (the case $a<b$ is similar). Then
\begin{align}\label{freeness}
(ab^{-1})\cdot z=z&\Longrightarrow z_{bn}^a=z_n=z_{bn}^b\quad\text{for all $n$}\\
&\Longrightarrow z_{bn}^{a-b}=1\quad\text{for all $n$}\notag\\
&\Longrightarrow z_{bn}=1\quad\text{for all $n$ \quad (because $a-b\in \nx$)}\notag\\
&\Longrightarrow z_n=z^b_{bn}=1\quad\text{for all $n$}.\notag
\end{align}
Thus the action is free on $\mathcal{S}\setminus\{1\}$, as claimed.
 
To establish minimality, we take $1\neq w\in\mathcal{S}$ and aim to prove that the orbit of $w$ is dense in $\mathcal{S}$.  To see this, we take $z\in\mathcal{S}$. Since the inverse limit is topologised as a subset of $\prod_{n\in \nx}\T$, and   a basic open neighbourhood $U$ of $z$ has the form $\big(\prod_{n\in F}U_n\big)\times\big(\prod_{n\notin F}\T\big)$ for some finite subset $F$ of $\nx$. We write $N:=\lcm\{n:n\in F\}$. Since each map $z_N\mapsto z_{N}^{n^{-1}N}$ is continuous for each $n\in F$, we can find a neighbourhood $V$ of $z_N$ such that $w_N\in V$ implies $w_n\in U_n$ for all $n\in F$. Then we replace $U$ by $V\times\big(\prod_{n\not=N}\T\big)$.

Now we need to find $a, b\in\nx$ such that $((ab^{-1})\cdot w)_N\in V$.  Since $w\neq 1$ there exists $b$ such that $w_b\neq 1$, and then $w_{bn}\not= 1$ for all $n\in \nx$. First, we suppose that there exists $w_{bN}\neq 1$ that is not a root of unity. Then $\{w_b^p:p\in\N\}$ is dense in $\T$, and we can choose $a\in \nx$ such that $w_{bN}^a\in N$. Then $((ab^{-1})\cdot w)_N=w_{bN}^a\in V$, as required.
 
The alternative is that all $w_{cN}\neq 1$ are roots of unity.  For $u\in\T$, we write $o(u)=n$ if and only if $u=e^{2\pi i k/n}$ for some $k$ co-prime to $n$. We claim that $\{o(w_{cN}):c\in\nx\}$ is unbounded. To see this, suppose not. Then there exists $a\in\nx$ such that $w^a_{cN}=1$ for all $c$. Then we have $(aN^{-1})\cdot w=1$ for all $c\in \nx$, which is impossible for $w\in \mathcal{S}\setminus\{1\}$ because the implication \eqref{freeness} implies that $w=1$.

Since $\{o(w_{cN})\}$ is unbounded, we can find $c\in \nx$ such that $\{z\in\T: z^{o(w_{cN})}=1\}\cap U\neq\emptyset$. Since
\[
\{z\in\T: z^{o(w_{cN})}=1\}=\{w_{cN}^a: a\in\N\},
\]
we get $a$ such that $w_{cN}^a\in U$.  Then $((a(cN)^{-1})\cdot w)_c=w_{cN}^a\in U$, as required. Thus the orbit of $w$ is dense in $\mathcal{S}$.
\end{proof}

\begin{prop}\label{primofC*(qxxq)}
We write $\phi$ for the homomorphism of $C^*(\Q)$ into $\C$ such that $\phi(u_r)=1$ for all $r\in \Q$. Then the primitive ideals of $\qx\ltimes_\alpha C^*(\Q)$ are $\{0\}$ and 
\[
\big\{\ker\big(\epsilon_\gamma\circ (\id\ltimes\phi)\big):\gamma\in \widehat{\qx}\big\}.
\]
The point $\{0\}$ is an open and dense subset of $\Prim\big(\qx\ltimes_\alpha C^*(\Q)\big)$, and the map $\gamma\mapsto \ker\big(\epsilon_\gamma\circ (\id\ltimes\phi)\big)$ is a homeomorphism of the compact space $\widehat{\qx}$ onto the complement of $\{0\}$. The ideal $\ker(\id\ltimes\phi)$ is simple.
\end{prop}

\begin{proof}
Lemma~\ref{free+min} says that the action of $\qx$ on $\mathcal{S}\setminus\{1\}$ is free and minimal, and hence it follows from \cite[Corollary to Theorem~2]{AS} that the crossed product $\qx\ltimes_{\alpha}C_0(\mathcal{S}\setminus\{1\})$ is simple. The exactness of \eqref{cpseqexact} implies that 
\[
\qx\ltimes_{\alpha}C_0(\mathcal{S}\setminus\{1\})=\ker(\id\ltimes\epsilon_1)=\ker(\id\ltimes\phi).
\]
The map $\gamma\mapsto \epsilon_\gamma$ is a homeomorphism of $\widehat{\qx}$ onto $\Prim C^*(\qx)$ (strictly speaking, the homeomorphism is $\gamma\mapsto \epsilon_\gamma\circ\F$). Thus exactness of \eqref{cpseqexact} implies that 
\[
\gamma\mapsto \ker\big(\epsilon_\gamma\circ(\id\ltimes\phi)\big)
\]
is a homeomorphism of $\widehat{ \qx}$ onto the closed subset
\[
\big\{P\in \Prim\big(\qx\ltimes_\alpha C^*(\Q)\big):\ker(id\ltimes\phi)\subset P\big\}
\]
(see \cite[Proposition~A.27]{tfb}, for example). The complement of this set in $\Prim\big(\qx\ltimes_\alpha C^*(\Q)\big)$ is
\[
\big\{P\cap \ker(id\ltimes\phi):\ker(id\ltimes\phi)\not\subset P\big\},
\]
which since $\ker(id\ltimes\phi)$ is simple is precisely $\{0\}$. This proves both that $\{0\}$ is primitive, and that the singleton set containing $\{0\}$ is open. Since the closure of $\{O\}$ is the set $\{P:O=\{0\}\subset P\}$ (see \cite[Definition~A.19]{tfb}),  it is all of $\Prim\big(\qx\ltimes_\alpha C^*(\Q)\big)$.
\end{proof}

\begin{proof}[Proof of Theorem~\ref{primcp}]
The homomorphism $\psi:\TT(\qp)\to \C$ is equivariant for the action $\alpha$ of $\qx$ on $\TT(\qp)$ and the trivial action of $\qx$ on $\C$. Thus there is a homomorphism $\id\ltimes\psi$, as asserted.

The exactness of the sequence \eqref{cpexact} implies that $P\mapsto (\id\ltimes\pi_u)^{-1}(P)$ is a homeomorphism of $\Prim\big(\qx\ltimes_\alpha C^*(\Q)\big)$ onto a closed subset of $\Prim\big(\qx\ltimes_\alpha \TT(\qp)\big)$. Since $I:=\ker(\id\ltimes \pi_u)=\qx\ltimes_{\alpha}\CC(\qp)$ is simple by Corollary~\ref{cpsimple}, the only ideal in $\qx\ltimes_\alpha \TT(\qp)$ that does not factor through $\id\ltimes\psi$ is $O=\{0\}$. It is the intersection of the primitive ideals containing it, and hence must itself be primitive. The set $\{O\}$ is open because the complement is closed, and dense because every primitive ideal contains $O$. The remaining assertions follow from Proposition~\ref{primofC*(qxxq)}.
 \end{proof}
 
It is now a relatively straightforward matter to use the isomorphism of Proposition~\ref{idcpToe} to convert  Theorem~\ref{primcp} into a parallel result about $\partial_{\mult}\TT(\nxxn)$.

\begin{theorem}\label{structuremultbdary}
There is a homomorphism $\pi$ of $\partial_{\mult}\TT(\nxxn)=C^*(\{T_{a,r}\})$ into $C^*(\qx)=C^*(\{u_a\})$ such that $\pi(T_{a,r})=u_a^*$. We denote by $q$ the quotient map from $\partial_{\mult}\TT(\nxnx)$ onto $\partial\TT(\nxnx)$. For $\gamma\in (\qx)^\wedge$, we  let $\epsilon_\gamma$ be the character of $C^*(\qx)$ such that $\epsilon_\gamma(u_a)=\gamma(a)$. 
\begin{enumerate}
\item The map $\gamma\mapsto \ker(\epsilon_\gamma\circ \pi)$ is a homeomorphism of the compact space $(\qx)^\wedge$ onto a closed subset $F$ of $\Prim\big(\partial_{\mult}\TT(\nxxn)\big)$; the ideal of $\partial_{\mult}\TT(\nxxn)$ with hull $F$ is $\ker\pi_t$. 
\item The other primitive ideals are $O:=\{0\}$ and $\ker q$. \item The set $\{O\}$ is open and dense in $\Prim\partial_{\mult}\TT(\nxxn)$, and the set $\{\ker q\}$ is open and dense in $\Prim\big(\partial_{\mult}\TT(\nxxn)\big)\setminus \{O\}$.
\item Both $\ker q$ and $\ker \pi/\ker q$ are simple.
\end{enumerate}
\end{theorem}

\begin{proof}
By Proposition~\ref{prop-commutingdiagram-revised}, the isomorphism $U\ltimes\pi_W: \qx\ltimes_\alpha \TT(\Q_+)\to \partial_{\mult}\TT(\nxxn)$  of Proposition~\ref{idcpToe}  carries the ideal $I:=\qx\ltimes_\alpha\CC(\qp)$ onto $\ker q$. 
The map $t:(a,r)\mapsto u_a^*$ is an isometric representation of $\nxxn$ in $C^*(\qx)$. Since each $u_a$ is unitary, $t$ satisfies $t_at_a^*=1$, and hence there is a homomorphism $\pi:=\pi_t$ on $\partial_{\mult}\TT(\nxxn)$, as claimed. The homomorphism $\id\ltimes\psi:  \qx\ltimes_\alpha \TT(\Q_+)\to C^*(\qx)$ of Theorem~\ref{primcp} corresponds to $\pi$, that is, $id\ltimes\psi=\pi\circ (U\ltimes\pi_W)$.  Thus $\ker\pi=\ker(\id\ltimes\psi)$.  Now the theorem follows from the parallel statements in Theorem~\ref{primcp}.
\end{proof}

\section{KMS states}

We consider the dynamics $\sigma:\R\to \Aut\TT(\nxxn)$ such that $\sigma_t(S)=S$ and $\sigma_t(V_a)=a^{it}V_a$ for $a\in \nx$. We want to study the KMS states of the system $(\TT(\nxxn),\sigma)$. We begin in this section by proving some results which are valid for all $\beta$. The first is very similar to what happens for $\TT(\nxnx)$: compare with \cite[Lemma~10.4]{LR-advmath}.

\begin{prop}
Every KMS$_\beta$ state of $(\TT(\nxxn),\sigma)$ factors through a state of the additive boundary quotient $\partial_{\add}\TT(\nxxn)$.
\end{prop}

\begin{proof}
Suppose that $\phi$ is a KMS$_\beta$ state of $(\TT(\nxxn),\sigma))$. Since the generator $S$ is fixed by the action $\alpha$, we have $1=\phi(S^*S)=\phi(SS^*)$. Thus $\phi$ vanishes on the element $1-SS^*$. Since $\alpha_t(1-SS^*)=1-SS^*$, we deduce from \cite[Lemma~6.2]{AaHR}, for example, that $\phi$ vanishes on the ideal generated by $1-SS^*$, and hence factors through a state of the quotient $\partial_{\add}\TT(\nxxn)$. 
\end{proof}

The next result is an analogue of \cite[Lemma~8.3]{LR-advmath}. As there, we write 
\[
S^{((k))}:=\begin{cases}S^k&\text{if $k\in \N$}\\
S^{*(-k)}&\text{if $k\in \Z$ and $k<0$.}
\end{cases}
\]
For our system $(\TT(\nxxn), \sigma)$, the spanning elements $V_aS^m(V_bS^n)^*$ satisfy \[\sigma_t(V_aS^m(V_bS^n)^*)=(ab^{-1})^{it} V_aS^m(V_bS^n)^*,\] and  hence are analytic.

\begin{prop}\label{charKMSprop}
Suppose that $\beta>0$ and $\phi$ is a state of $\TT(\nxxn)$. Then $\phi$ is a KMS$_\beta$ state of the system $(\TT(\nxxn),\sigma)$ if and only if 
\begin{equation}\label{charKMS}
\phi\big(V_aS^m(V_bS^n)^*\big)=\delta_{a,b}a^{-\beta}\phi\big(S^{((m-n))}\big)\quad\text{for all $(a,m)$, $(b,n)$ in $\nxxn$.}
\end{equation}
\end{prop} 

The following result is essentially proved in the proof of \cite[Lemma~8.3]{LR-advmath} (look in the second last paragraph of that proof).

\begin{lemma}\label{a'=b'}
Suppose that $a,b,c,d\in \Z$ satisfy $ac=bd$. Write 
\[
a=a'\gcd(a,d),\quad b=b'\gcd(b,c),\quad c=c'\gcd(b,c)\quad\text{and}\quad d=d'\gcd(a,d).
\]
Then $a'=b'$ and $c'=d'$. 
\end{lemma}

\begin{proof}[Proof of Proposition~\ref{charKMSprop}]
Suppose that $\pi$ is a KMS$_\beta$ state of $(\TT(\nxxn),\alpha)$ and $V_aS^m(V_bS^n)^*$ is a fixed spanning element. Then the KMS condition gives
\begin{align}\label{useKMS1}
\phi\big(V_aS^m(V_bS^n)^*\big)&=\phi\big((V_bS^n)^*\alpha_{i\beta}(V_aS^m)\big)\\
&=a^{-\beta}\phi\big((V_bS^n)^*V_aS^m\big)\notag\\
&=a^{-\beta}b^{\beta}\phi\big(V_aS^m(V_bS^n)^*\big)\notag.
\end{align}
Since $\beta>0$ and $V_a^*V_a=1$, we have
\begin{equation}\label{useKMS2}
\phi\big(V_aS^m(V_bS^n)^*\big)\not=0\Longrightarrow b^{\beta}=a^{\beta}\Longrightarrow b=a.
\end{equation}
Together, \eqref{useKMS1} and \eqref{useKMS2} imply \eqref{charKMS}.

Now we suppose that $\phi$ satisfies \eqref{charKMS}, and aim to prove that $\phi$ is a KMS$_\beta$ state. Since the elements $V_aS^m(V_bS^n)^*$ are analytic and span a dense subalgera, it suffices to prove that 
\[
\phi(xy)=\phi(y\alpha_{i\beta}(x))
\]
for $x=V_aS^m(V_bS^n)^*$ and $y=V_cS^p(V_dS^q)^*$. 
We write $a=a'\gcd(a,d),\quad b=b'\gcd(b,c),\quad c=c'\gcd(b,c)\quad\text{and}\quad d=d'\gcd(a,d)$  as in Lemma~\ref{a'=b'} above. 
On one hand, we have 
\begin{align*}
xy&=V_aS^mS^{n*}V_b^*V_cS^pS^{q*}V_d^*\\
&=V_aS^mS^{n*}(V_{b'}^*V_{c'})S^pS^{q*}V_d^*\quad\text{since $V_{\gcd(b,c)}^*V_{\gcd(b,c)}=1$}\\
&=V_aS^mS^{n*}(V_{c'}V_{b'}^*)S^pS^{q*}V_d^*\quad\text{using (T3)}\\
&=(V_aS^mV_{c'}S^{c'n*})(S^{b'p}V_{b'}^*S^{q*}V_d^*)\quad\text{using (T4)}\\
&=(V_{ac'}S^{c'm}S^{c'n*})(S^{b'p}S^{b'q*}V_{b'd})\quad\text{using (T1).}
\end{align*}
Thus \eqref{charKMS} implies that 
\begin{equation}
\label{compLHS}\phi(xy)=\delta_{ac',b'd}(ac')^{-\beta}\phi\big(S^{((c'm+b'p-c'n-b'q))}\big).
\end{equation}
On the other hand, we have
\begin{equation}\label{compRHS}
\phi(yx)=\delta_{ca',bd'}(ca')^{-\beta}\phi\big(S^{((a'p+d'm-a'q-d'n))}\big),
\end{equation}
Next, we observe that
\begin{align*}
ca'=bd'&\Longleftrightarrow ca=\gcd(a,d)ca'=\gcd(a,d)bd'=bd\\
&\Longleftrightarrow
ac'=\gcd(b,c)^{-1}ca=\gcd(b,c)^{-1}bd=db',
\end{align*}
so the Kronecker delta functions in \eqref{compLHS} and \eqref{compRHS} are the same. Thus $\phi(xy)=0=\phi(xy)$ when $ca'\not= bd'$. When $ca'= bd'$, we have 
\[
\phi(y\alpha_{i\beta}(x))=a^{-\beta}b^\beta \phi(yx)=a^{-\beta}b^\beta(ca')^{-\beta}\phi\big(S^{((a'p+d'm-a'q-d'n))}\big).
\]
By Lemma~\ref{a'=b'} we have $a'=b'$. Therefore the coefficient on the right-hand side is
\[
\Big(\frac{aca'}{b}\Big)^{\-\beta}=\Big(\frac{a\gcd{b,c}a'}{b}\Big)^{\-\beta}=\Big(\frac{ac'a'}{b'}\Big)^{\-\beta}=(ac')^{\beta}
\]
which is the same as the coefficient in the formula \eqref{compLHS}. 
Finally we have (using again that $a'=b'$ and $c'=d'$)
\begin{align*}
a'p+d'm-a'q-d'n&=a'(p-q)+d'(m-n)=b'(p-q)+c'(m-n)\\
&=c'm+b'p-c'n-b'q.
\end{align*}
We deduce that when $ac'=b'd$ we have
\[
\phi\big(y\alpha_{i\beta}(x)\big)=(ac')^{-\beta}\phi\big(S^{((c'm+b'p-c'n-b'q))}\big)=\phi(xy), 
\]
as required.
\end{proof}

\section{KMS states for large inverse temperatures}\label{sec-KMS-largebeta}

We consider now KMS$_\beta$ states of $(\TT(\nxxn),\sigma)$ for $\beta>1$. We begin by stating our main result.

\begin{theorem}\label{mainKMSthm}
Let $\sigma:\R\to \Aut\TT(\nxxn)$  be the dynamics on $\TT(\nxxn)$ determined by $\sigma_t(S)=S$ and $\sigma_t(V_a)=a^{it}V_a$ for $t\in \R$ and $a\in \nx$.
\begin{enumerate}
\item
For each state $\varphi$ on the Toeplitz algebra $\TT = C^*(S)$ there is a unique ground 
state $\omega_\varphi $ of $\TT(\nxxn)$ such that 
\begin{equation}\label{omegamudef}
\omega_\varphi( V_a S^mS^{*n} V_b^*) =  \delta_{a,b} \delta_{a,1} \varphi(S^mS^{*n})\quad\text{for $a,b\in \nx$ and $m,n\in\N$}.
\end{equation} 
The map $\varphi\mapsto \omega_\varphi$ 
is an affine w*-homeomorphism of the set of 
states of $\TT$  onto the ground states of $(\TT(\nxxn),\sigma)$.

\item For each probability measure $\mu$ on the circle and each $\beta >1$ there is a unique KMS$_\beta$ state of $\TT(\nxxn)$ such that 
\begin{equation}\label{kmstaubetaformula}
\psi_{\mu,\beta} ( V_a S^mS^{*n} V_b^*   ) = \delta_{a,b} \frac{a^{-\beta}}{\zeta(\beta)} \sum_{c\in \nx} c^{-\beta}\int_\T z^{c(m-n)} d\mu(z).
\end{equation}
The map $\mu\mapsto \psi_{\mu,\beta}$ is an affine w*-homeomorphism of the simplex of  probability measures on $\T$
onto the simplex of KMS$_\beta$ states of $(\TT(\nxxn),\sigma)$.
\end{enumerate}
\end{theorem}

Part (2) of this theorem is very similar to part (3) of \cite[Theorem~7.1]{LR-advmath}, which is about KMS$_\beta$ states of $\TT(\nxnx)$ for $\beta\in (2,\infty]$. Our proof relies on general results from \cite{diri} about KMS states on semigroup crossed products. 
So we begin by writing $\TT(\nxxn)$ as a semigroup crossed product. The next proposition describes the underlying algebra $\mathfrak C$ of this crossed product, and the following one describes the crossed-product decomposition.

\begin{prop}\label{rightToeplitz-structure}
There is an action $\theta$ of $\widehat \qps$ by automorphisms of $\TT(\nxxn)$ such that 
\begin{equation*}
\theta_\chi(S)=S\quad\text{and}\quad \theta_\chi (V_a) = \chi(a) V_a\text{ for $\chi \in \widehat \qps$ and $a\in \nx$.}
\end{equation*}  
The fixed-point algebra of this action is 
\begin{equation}\label{spanning2}
\mathfrak C:=\TT(\nxxn)^\theta =  \clsp\{ V_a S^m S^{*n} V_a^*: a\in \nx ,  m, n \in \N\},
\end{equation}
and the projection $V_aV_a^*$ is central in $\mathfrak C$ for each $a\in \nx$.
\end{prop}

\begin{proof}
If $s$ and $\{v_a: a\in \nx\}$ are universal generators for $\TT(\nxxn)$, then so are $s$ and $\{\chi(a) v_a: a\in \nx\}$ for each $\chi \in \widehat \qps$, and hence the universal property gives an automorphism $\theta_\chi$. The universal property also implies that $\chi\mapsto \theta_\chi$ is a group homomorphism, and because the automorphisms are all isometric, it suffices to check continuity on the spanning family \eqref{spanning}, and this is straightforward. The fixed-point algebra of $\theta$ is the range of the associated conditional expectation $E^\theta$ on $\TT(\nxxn)$, and hence \eqref{spanning2} follows from \eqref{spanning}.

Now we suppose that $a\in \nx$, and aim to prove that $V_aV_a^*$ is central in $\TT(\nxxn)^\theta$. Let  $b \in \nx$ and $m,n \in \N$. We write $a = a' \gcd(a,b)$, $b = b'\gcd(a,b)$ with $\gcd(a',b') =1$ and  $ ab' = ba'= \lcm(a,b)$. Using first (T3) and then the adjoints of (T1) and (T4), we find that
\begin{align}\label{compute a'b'}
V_aV_a^* (V_b S^m S^{*n}V_b^*)&= V_a V_{b'} V_{a'}^*S^m S^{*n}V_b^* = V_{ab'} S^{a'm} S^{*a'n}V_{a'}^* V_b^*\\
&= V_{\lcm(a,b)} S^{a'm} S^{*a'n}V_{\lcm(a,b)}^*.\notag
\end{align}
Similarly, using (T1) and (T4) directly gives 
\[
 (V_b S^m S^{*n}V_b^*) V_aV_a^* =V_bS^m S^{*n}V_{a'}V_{b'}^* V_a^*= V_{ba'} S^{a'm} S^{*a'n}V_{ab'}^*= V_{\lcm(a,b)} S^{a'm} S^{*a'n}V_{\lcm(a,b)}^*.
 \]
Thus 
\[
V_aV_a^* (V_b S^m S^{*n}V_b^*)= V_{\lcm(a,b)} S^{a'm} S^{*a'n}V_{\lcm(a,b)}^*= (V_b S^m S^{*n}V_b^*) V_aV_a^*,
\] 
and $V_a V_a^*$ is central as claimed.
\end{proof}

The next proposition describes $\TT(\nxxn)$ as a semigroup crossed product to which Theorem~12 of \cite{diri} applies.

\begin{proposition}\label{Ctimesnx}
Let  $\mathfrak C $  be the subalgebra of $\TT(\nxxn)$ defined in \eqref{spanning2}. For $a\in \nx$ , we define $\alpha_a:\mathfrak C\to \mathfrak C$ by $\alpha_a(X)=V_a X V_a^*$. 
\begin{enumerate}
\item Then $\{\alpha_a: a\in \nx\}$ is a semigroup of endomorphisms of $\mathfrak C$, we have  $\alpha_a (\mathfrak C) = V_aV_a^* \mathfrak C V_aV_a^*$,  and $\alpha_a$ has a left inverse given by $\theta(Y)=V_a^* Y V_a$.
\item There is a canonical isomorphism
 \begin{equation}\label{TTisc-p}
\TT(\nxxn) \cong \mathfrak C \rtimes_\alpha \nx
\end{equation}
which matches up the two copies of $\mathfrak C$, and for $a\in\nx$  takes $V_a \in\TT(\nxxn)$ to the generating isometry $w_a\in  \mathfrak C \rtimes_\alpha \nx$. 
\item The semigroup action respects the lattice structure, as defined in \cite[Definition 3]{diri}.
\end{enumerate}
\end{proposition}

\begin{proof} For $a,b\in \nx $ and $m,n \in \Z$ we have 
\[
\alpha_a (V_b S^m S^{*n}V_b^*) = V_{a}(V_bS^m S^{*n}V_b^*) V_{ a}^*=V_aV_a^*(V_{ab} S^m S^{*n}V_{ab}^*)V_aV_a^*.
\]
Thus  $\alpha_a$ maps $\mathfrak C$ to  $\mathfrak C$ and  $\alpha_a (\mathfrak C)\subset V_aV_a^* \mathfrak C V_aV_a^*$. It follows from \eqref{compute a'b'} that $\alpha_a (\mathfrak C) = V_aV_a^* \mathfrak C V_aV_a^*$. 
Since $V$ is a semigroup of isometries, $\{\alpha_a\}$ is a semigroup of  endomorphisms of $\mathfrak C$ with left inverses as described.

To prove the existence of the isomorphism \eqref{TTisc-p}, we observe that the inclusion $ \iota: \mathfrak C \subset \TT(\nxxn)$ 
and  the semigroup of isometries $\{V_a: a\in \nx\}$ form a covariant pair for the semigroup dynamical system 
$(\mathfrak C, \nx, \alpha)$. Thus we  have a homomorphism  $ \iota \rtimes V : \mathfrak C \rtimes_\alpha \nx \to \TT(\nxxn)$. 
This proves in particular that the image of $\mathfrak C$ in $\mathfrak C \rtimes_\alpha \nx$ is faithful and that $w_a$ maps to $V_a$. 

To obtain the inverse of $\iota \rtimes V$, we show that the relations (T0)--(T4) hold for the isometry $S\in \mathfrak C$ and the canonical semigroup of isometries $w_a$ viewed as elements of  $\mathfrak C \rtimes_\alpha \nx$. 
Relations (T0) and (T2) reflect that $S$ and $w$ are isometric representations. For  
 (T3), let $a,b\in\nx$. First notice that covariance implies that   
$w_aw_a^*=\alpha_a(1)=V_aV_a^*.$  Thus
\begin{align*}
w_aw_a^* w_bw_b^* &=\alpha_a(1) \alpha_b(1) = V_aV_a^* V_bV_b^* \\
&= V_{\lcm(a,b)}V^*_{\lcm(a,b)}\quad \text{by  \eqref{compute a'b'} with $m=n=0$}
\\&= \alpha_{\lcm(a,b)} (1) = w_{\lcm(a,b)}w_{\lcm(a,b)}^* .
\end{align*}
Now suppose that $\gcd(a,b)=1$, so that $\lcm(a,b)=ab$. Then 
\[
w_a^*w_b=w_a^*(w_aw_a^* w_bw_b^*)w_b=w_a^*w_{\lcm(a,b)}w_{\lcm(a,b)}^*w_b=w_bw_a^*,
\] 
giving (T3).

To prove that (T1) holds in the crossed product, 
we recall that $V_a^* S V_a = S^a$ in the C*-algebra $\mathfrak C$ because (T4) holds in $\TT(\nxxn)$. Thus 
\begin{align*}
w_a^* S w_a&= w_a^* (w_aw_a^* S w_aw_a^*) w_a = 
w_a^* (V_aV_a^*S V_aV_a^*) w_a \\ 
&=w_a^*  V_a S^aV_a^* w_a = w_a^*  \alpha_a (S^a) w_a 
= w_a^* w_a S^a w_a^* w_a = S^a.
\end{align*}
Now
\begin{align*}
Sw_a= S w_a w_a^* w_a &= S \alpha_a(1) w_a = S(V_aV_a^*)w_a\\
&=V_aS^aV_a^*w_a\quad\text{using (T1) for $(S,V)$}\\
&=V_aV_a^*Sw_a\quad\text{using (T4) for $(S,V)$}\\
&=\alpha_a(1) S w_a =w_a w_a^* S w_a \\
&= w_a S^a,
\end{align*}
giving (T1). A similar argument with $S^*$ in place of $S$ shows that (T4) holds in the crossed product. Using the universal property of $\TT(\nxxn)$ from \proref{presentation} shows that the
assignment $S \mapsto S$ and $V_a\mapsto w_a$ gives an inverse for $\iota \rtimes V$.

For  $\alpha$  to respect the lattice structure in the sense of \cite[Definition~3]{diri} requires that $\alpha_a(\mathfrak C)$ is an ideal in $\mathfrak C$ and that $\alpha_a(1)\alpha_b(1)=\alpha_{a\vee b}(1)$ for $a,b\in\nx$.   The first follows because $\alpha_a(\mathfrak C)=V_a V_a^*\mathfrak CV_aV_a^*$ and $V_aV_a^*$ is central in $\mathfrak C$.  The second follows because 
\[
\alpha_a(1)\alpha_b(1)=V_aV_a^*V_bV_b^*=V_{\lcm(a,b)}V ^*_{\lcm(a,b)}=\alpha_{\lcm(a,b)}(1),
\]
and because $a\vee b=\lcm(a,b)$ in the lattice-ordered semigroup $\nx$.
 \end{proof}
 
\begin{proof}[Proof of Theorem~\ref{mainKMSthm}]
For each $p$ in the set $\primes$ of primes, denote by $\nx_p$ the set of numbers with all prime factors less than $p$. We now define 
\[
\mathfrak C_{ p} := \clsp\{ V_a S^mS^{*n} V_a^*: a \in \nx_p, m,n\in \N\}.
\] 
A computation like that in \eqref{compute a'b'} gives 
 \begin{equation}\label{Cpalg}
 (V_a S^mS^{*n} V_a^*)(V_b S^jS^{*k} V_b^*) =  V_a S^mS^{*n} (V_{b'} V_{a'}^*)S^jS^{*k} V_b^* =  V_{ab'} S^{mb'}S^{*nb'}  S^{ja'}S^{*ka'} V_{ba'}^*;
 \end{equation}
since $S^{*c}S^d$ is either $S^{d-c}$ or $S^{*(c-d)}$, and since the prime factors of $ab' = ba' = \lcm(a,b)$ appear in $a$ or $b$,  we deduce from \eqref{Cpalg} that 
$\mathfrak C_p$ is a C*-subalgebra of $\mathfrak C$. Since $\nx_q\subset \nx_p$ for $q\leq p$, we have $\mathfrak C_q \subset \mathfrak C_p$  for $q\leq p$, and 
\[
\bigcup_{p\in \primes} \mathfrak C_p =\lsp\{V_aS^mS^{*n}V_a^*:a\in\nx,\text{ and }m,n\in \N\big\}
\]
is dense in $\mathfrak C$.  For $p\in \primes$, we consider the projection 
\[
Q_p := \prod_{ q\leq p, q\in \primes} (1 - V_{q}V_{q}^*); 
\] 
since Proposition~\ref{rightToeplitz-structure} says that each $V_q V_q^*$ is central in $\mathfrak C$, the projection  $Q_p$ is also central  in $\mathfrak C_p$.
Let $a\in \nx_p \setminus \{1\}$. Then $a$ has at least one prime factor $q \leq p$. Since $V$ is Nica-covariant for $\nx$, we have
\[
V_qV_q^*V_a=V_qV_q^*V_aV_a^*V_a=V_{q\vee a}V^*_{q\vee a}V_a=V_aV_a^*V_a=V_a.
\]
Then
\[
Q_p V_a S^mS^{*n}V_a^* = \big(Q_p(1-V_qV_q^*)\big)V_a S^mS^{*n}V_a^*=\big(Q_p(1-V_qV_q^*)\big)\big(V_qV_q^* V_a\big) S^mS^{*n}V_a^*= 0.
\]
Hence we have  
\[
Q_p \mathfrak C_p Q_p=\clsp\{Q_p S^mS^{*n}Q_p \}.
\]
 
The relations (T1) and (T4) imply that $SV_qV_q^*=V_qV_q^*S$. Thus $Q_p S Q_p=Q_pS=SQ_p$.
We claim that the map $S \mapsto Q_p S Q_p$ extends to an isomorphism 
$\pi_{Q_pS}$ of $\TT = C^*(S)$ onto $Q_p \mathfrak C_p Q_p$. To see this, we observe that
\[
(Q_pS)^* (SQ_p) = (Q_pS^*Q_p) (Q_pSQ_p) = (Q_pS^*)(SQ_p) = Q_p,
\]
and hence $Q_pS=Q_pSQ_p$ is an isometry in the corner $Q_p\mathfrak C_p Q_p$. On the other hand, since $Q_p\varepsilon_{1,0}=\varepsilon_{1,0}$ and $(Q_pS)^*\varepsilon_{1,0}=0$,  $Q_pS$ is a proper isometry in $\mathfrak C_p$. Thus Coburn's theorem gives the required isomorphism $\pi_{Q_pS}$. In the previous paragraph we saw that $Q_pS$ generates $Q_p\mathfrak C_p Q_p$, and hence the range of $\pi_{Q_pS}$ is all of $Q_p \mathfrak C_p Q_p$.

Let  $\varphi$ be a state of $\TT$.  For each prime $p$, we get a state $\varphi\circ\pi_{Q_pS}^{-1}$ on $Q_p\mathfrak C_p Q_p$ such that $\varphi\circ\pi_{Q_pS}^{-1}(Q_p)=\phi(1)=1$. Now we extend $\varphi\circ\pi_{Q_pS}^{-1}$ to  a positive functional $\omega_p$ on  $ \mathfrak C_p$ by setting  
\[
\omega_p(c) =\varphi\circ\pi_{Q_pS}^{-1}(Q_p c Q_p)\text{\quad for $c\in \mathfrak C_p$.}
\]
Since $\omega_p(1)=1$,  $\omega_p$ is a state (by \cite[II.6.2.5]{B}, for example). Then for $a\in \N_p^\times$ we have
\begin{align}
\label{formomegap}\omega_p(V_a S^mS^{*n} V_a^*)&=\varphi\circ\pi_{Q_pS}^{-1}\big(Q_p(V_a S^mS^{*n} V_a^*)Q_p\big)\\
&=\begin{cases}\varphi\circ\pi_{Q_pS}^{-1}(0)\text{\quad if $a\not=1$, and}\\
\varphi\circ\pi_{Q_pS}^{-1}(Q_pS^mS^{*n}Q_p)\text{\quad if $a=1$}
\end{cases}\notag\\
&=\delta_{a,1} \varphi(S^m S^{*n}).\notag
\end{align}
For $q\leq p$, we have $\N_q^\times\subset \N_p^\times$, and hence \eqref{formomegap} implies that for $a\in \N_q^\times$ we have
\[
\omega_p(V_a S^mS^{*n} V_a^*)=\delta_{a,1}(\varphi(S^mS^{*n})=\omega_q(V_a S^mS^{*n} V_a^*).
\]
Thus $\omega_p |_{\mathfrak C_q} = \omega_q$, and we can define $\omega_\infty:\bigcup_{p}\mathfrak C_p\to \C$ by $\omega_\infty(c)=\omega_p(c)$ for $c\in \mathfrak C_p$. Since each $\omega_p$ is a linear functional with norm $1$, $\omega_\infty$ also has norm $1$, and extends to a linear functional on the closure $\mathfrak C$; since each $\omega_p$ is positive, so is $\omega_\infty$. Since $1\in \mathfrak C$ and $\omega_p(1)=$ for all $p$, we have $\omega_\infty(1)=1$. Thus $\omega_\infty$ is a state.

We now consider the action $\theta$ of Proposition~\ref{rightToeplitz-structure}, and the associated expectation $E^\theta$ of $\TT(\nxxn)$ onto $\TT(\nxxn)^\theta$.
Then the state $\omega_\varphi:= \omega_\infty \circ E^\theta$ of $\TT(\nxxn)$ satisfies \eqref{omegamudef}. It follows from \eqref{omegamudef} that the map $\varphi \mapsto \omega_\varphi \circ E^\theta$ is injective, affine and weak*-continuous. 

To complete the proof of part~(1), it remains to show that $\omega_\varphi$ is a ground state for every $\varphi$ and that every ground state arises this way.
From \proref{Ctimesnx} we know that the semigroup $\{\alpha_a\}$ respects the lattice structure, and that the dynamics on $\TT(\nxxn) \cong \mathfrak C \rtimes_\alpha \nx $ has the form $\sigma^N$ of \cite[Theorem 12]{diri} for the function $N:\nx\to (0,\infty)$ defined by
 $N(a) = a$ for $a\in \nx$. Thus the hypothesis of  \cite[Theorem 12]{diri} are satisfied, and part (2) of that theorem says that a state $\omega$ of $(\TT(\nxxn), \sigma)$ is a ground state if and only if it factors through the conditional expectation $E^\theta$ and its restriction to the fixed point algebra $\mathfrak C$ vanishes on all the corners $\alpha_a(\mathfrak C)$ for $a\neq 1$, or equivalently
if and only if
\begin{equation*}
\omega( V_a S^mS^{*n} V_b^*) = \delta_{a,b}\delta_{a,1} \omega(S^mS^{*n})\quad\text{for all $a,b\in \nx$ and $m,n\in \N$},
\end{equation*}
which is \eqref{omegamudef}.
Thus $\omega$ is a ground state if and only if it 
is equal to $\omega_\varphi$ for $\varphi := \omega |_{\TT}$. We have now proved part~(1).

\smallskip

For part (2), fix $\beta>1$. We seek to apply \cite[Theorem 20]{diri} to $\TT(\nxxn) \cong \mathfrak C \rtimes_\alpha \nx$ with  function $N(a)=a$. The critical inverse temperature is $1$. As we noted above, the hypotheses of  \cite[Theorem~12]{diri} are satisfied, and these are also hypotheses of  \cite[Theorem~20]{diri}. We also need to know that the isometries $\{ V_a : a\in \nx\}$ generate a faithful copy of $\TT(\nx)$ in $\TT(\nxxn)$ and that their ranges $V_aV_a^* = \alpha_a(1)$ generate a faithful copy of the diagonal subalgebra $B_{\nx}$; this follows from \cite[Theorem~3.7]{quasilat}.
Thus  \cite[Theorem 20]{diri}  gives an affine weak* homeomorphism $T_\beta$ from the space of tracial ground states  of $\mathfrak C \rtimes_\alpha \nx \cong \TT(\nxxn)$ onto the  KMS$_\beta$ states of $\mathfrak C \rtimes_\alpha \nx$. The normalising factor $\zeta_N(\beta)$ in that theorem is here the Riemann zeta function.  Thus the formula for $T_\beta$ is 
\[
T_\beta \omega =\frac{1}{\zeta(\beta)}\sum_{d\in\nx}d^{-\beta} \omega\circ\gamma_d,
\]
where $\gamma$ is the left inverse of  $\alpha$. 

By  part (1), the tracial ground states of $\TT(\nxxn)$ are parametrised by traces on the Toeplitz algebra $\TT=C^*(S)$.  If $\phi$ is a trace on $\TT$, then $\phi$ factors through a state of of the quotient $C(\T)$ by the relation $1-SS^*=0$, and hence $\phi$ is given  by a probability measure on the circle\footnote{To see this, note  (following the argument of \cite[Lemma~2.2]{aHLRS1}) that the tracial property gives $\phi(SS^*)=\phi(S^*S)=\phi(1)=1$. Thus with $p:=1-SS^*$, we have $\phi(p)=0$. Now for any $a\in \TT$ we have
\[
0\leq\phi(pa^*ap)\leq \phi(p\|a\|^2p)=\|a\|^2\phi(p)=0,
\]
and since the positive elements span the ideal generated by $p$, it follows that $\phi$ vanishes on $p\TT p$. But now for any $a,b\in \TT$, we have
\[
\phi(apb)=\phi((ap)(pb))=\phi((pb)(ap))=0,
\]
and $\phi$ vanishes on the ideal generated by $p$.}. 

It remains to establish  \eqref{kmstaubetaformula}.   The formula for $T_\beta\omega$, when applied to the tracial ground state $\omega = \omega_{\tau_\mu}$ arising from a probability measure $\mu$ on $\T$, becomes
\begin{equation}\label{kmsformula}
\psi_{\tau_\mu,\beta} (V_a S^mS^{*n} V_b^*) = \frac{1}{\zeta(\beta)} \sum_{d \in \nx} d^{-\beta} \omega_{\tau_\mu}  (V_d^* V_a S^mS^{*n} V_b^* V_d).
\end{equation}
Write 
$a = a'\gcd(a,d)$, 
$d = d'\gcd(a,d)$ and $b=b'\gcd(d,b)$.
Then 
\begin{align*}
V_d^*V_aS^mS^{*n} V_b^*V_d&=V_{d'}^*V_{a'}S^mS^{*n} V_{b'}*V_{d'}=V_{a'}V_{d'}^*S^mS^{*n} V_{d'}V_{b'}^* \quad\text{using (T3)}\\
&=V_{a'}S^{md'}V_{d'}^*V_{d'}S^{*nd'} V_{d'}V_{b'}* =V_{a'}S^{md'}S^{*nd'} V_{d'}V_{b'}^* \quad\text{using (T4)}.
\end{align*} 
But $\omega_{\tau_\mu} (V_d^*V_aS^mS^{*n} V_b^*V_d)= \delta_{a,b}\delta_{a',1}  \omega_{\tau_\mu} (S^mS^{*n})$ by \eqref{omegamudef}. 
So  the only nonzero terms on the right-hand side of \eqref{kmsformula} are the ones where $a'=1$, equivalently the ones where $d$ is  a multiple of $a$. Thus
\begin{align*} 
\psi_{\tau_\mu,\beta} (V_a S^mS^{*n} V_b^*) &= \delta_{a,b}\frac{1}{\zeta(\beta)} \sum_{d \in \nx, a\mid d} d^{-\beta} \omega_{\tau_\mu}  (V_d^* V_a S^mS^{*n} V_a^* V_d)\\
&=
 \delta_{a,b} \frac{1}{\zeta(\beta)} \sum_{c \in \nx} (ac)^{-\beta}\omega_{\tau_\mu}( V_{c}^*S^mS^{*n}V_{c})\\
&=
 \delta_{a,b} \frac{a^{-\beta}}{\zeta(\beta)} \sum_{c \in \nx} c^{-\beta}w_{\tau_\mu}( S^{cm}S^{*cn})\\
  &=   \delta_{a,b}\frac{a^{-\beta}}{\zeta(\beta)} \sum_{c \in \nx} c^{-\beta}{\tau_\mu}( S^{cm}S^{*cn})\\
    &=  \delta_{a,b} \frac{a^{-\beta}}{\zeta(\beta)} \sum_{c \in \nx} c^{-\beta}\int_\T z^{c(m-n)} d\mu(z).\qedhere
\end{align*}
\end{proof}

\begin{remark}
This proof is quite different from the one for $\nxnx$ in \cite{LR-advmath}, which uses  a Hilbert-space representation. We can also realise the KMS$_\beta$ states of $\TT(\nxxn)$ spatially, and we outline the construction. For a probability measure $\mu$ on $\T$, we consider the Hilbert space $L^2(\T,\mu)$, and define $H_\mu:=\ell^2(\nx,L^2(\T,\mu))$. For $f\in L^2(\T,\mu)$ and $c\in \nx$, we define $fe_d:c\mapsto \delta_{c,d}f$. Then the elements $\{fe_d\}$ span a dense subspace of $H_\mu$. We write $M$ for the representation of $C(\T)$ on $L^2(\T,\mu)$ by multiplication operators. Then the formulas $S(fe_d)=(M^df)e_c$ and $V_a(fe_d)=fe_{ad}$ define isometries on $H_\mu$ which satisfy the relations (T1)--(T4) of Proposition~\ref{presentation}, and hence give a representation $\pi_\mu:\TT(\nxxn)\to B(H_\mu)$. Then
\[
\psi_{\beta,\mu}(T)=\zeta(\beta)^{-1}\sum_{d\in \nx}d^{-\beta}\big(\pi_\mu(T)(1e_d)\,|\,1e_d\big)
\]
gives a specific formula for a state, and we can deduce from Proposition~\ref{charKMSprop} that it is a KMS state. The details are messy, but easier than the ones in \cite[\S9]{LR-advmath}, and we wind up with the formula \eqref{kmstaubetaformula}, and hence we are describing the same KMS states.

However, we were not able to adapt the arguments of \cite[\S10]{LR-advmath} to prove that every KMS$_\beta$ state has the form $\psi_{\beta,\mu}$. So the extra information that we get from using the results in \cite{diri} is the surjectivity of the parametrisation in \cite[Theorem~20]{diri}.
\end{remark}

\section{KMS states with inverse temperature $1$}

For each probability measure $\mu$ on $\T$, compactness of the state space implies that there is a sequence $\beta_n\to 1+$ such that $\{\psi_{\beta_n,\mu}\}$ converges weak* to a state $\psi_{1,\mu}$. It follows on general grounds \cite[Proposition~5.3.23]{bra-rob} that the limit $\psi_{1,\mu}$ is a KMS$_1$ state of $(\TT(\nxxn),\sigma)$ (or one can apply Proposition~\ref{charKMSprop}). We naturally ask whether it is the only KMS$_1$ state, as in \cite{LR-advmath} for the usual left Toeplitz algebra $\TT(\nxnx)$. 

We answer this by looking at some specific measures where we can compute the values of the
state explicitly. We leave the detailed study of critical and subcritical equilibria to future work.

\begin{example}
We take for $\mu$ the normalised Lebesgue measure on $\T$, so that
\[
\int_{\T} f\,d\mu=\int^1_0 f(e^{2\pi ix})\,dx.
\]
Then we have
\[
\int_{\T} z^{c(m-n)}\,d\mu(z)=\begin{cases}
1&\text{if $c(m-n)=0$, and}\\
0&\text{if $c(m-n)\not=0$.}
\end{cases}
\]
So we have
\[ 
\psi_{\beta,\mu}\big(V_aS^mS^{*n}V_b^*\big)=\delta_{a,b}\delta_{m,n}\zeta(\beta)^{-1}a^{-\beta}\sum_{c\in\nx}c^{-\beta}=\delta_{a,b}\delta_{m,n}a^{-\beta}.
\]
In the limit as $\beta\to 1^+$, we find
\[ 
\psi_{1,\mu}\big(V_aS^mS^{*n}V_b^*\big)=\delta_{a,b}\delta_{m,n}a^{-1}.
\]
\end{example}

\begin{example}
For $\mu=\delta_1$ the point mass at $1\in \T$, we have
\[
\int_{\T} z^{(m-n)c}\,d\delta_1(z)=1^{c(m-n)}=1\quad\text{for all $m,n$ and $c$.}
\]
Thus we have
\[
\psi_{\beta,\delta_1}\big(V_aS^mS^{*n}V_b^*\big)=\delta_{a,b}\zeta(\beta)^{-1}a^{-\beta}\sum_{c\in\nx}c^{-\beta}=\delta_{a,b}a^{-\beta}\quad\text{for all $m,n\in\N$,}
\]
and 
\[
\psi_{1,\delta_1}\big(V_aS^mS^{*n}V_b^*\big)=\delta_{a,b}a^{-1}\quad\text{for all $m,n\in\N$.}
\]
\end{example}

\begin{example}
For $\mu=\delta_{-1}$, we have
\[
\int_{\T} z^{c(m-n)}\,d\delta_{-1}(z)=(-1)^{c(m-n)}\quad\text{for all $m,n$ and $c$.}
\]
Thus we have
\[
\psi_{\beta,\delta_{-1}}\big(V_aS^mS^{*n}V_b^*\big)=\delta_{a,b}\zeta(\beta)^{-1}a^{-\beta}\sum_{c\in\nx}c^{-\beta}(-1)^{c(m-n)}.
\]
We now take $a=b$, so that the delta function disappears. If $m-n$ is even, then as in the previous example,
\[
\psi_{\beta,\delta_{-1}}\big(V_aS^mS^{*n}V_a^*\big)=\zeta(\beta)^{-1}a^{-\beta}\sum_{d\in\nx}c^{-\beta}=a^{-\beta}.
\]
If $m-n$ is odd, then we note that 
\[
\sum_{c\in\nx\setminus 2\nx}c^{-\beta}=\sum_{c\in\nx}c^{-\beta}-\sum_{c\in\nx}(2c)^{-\beta}
\]
and hence 
\begin{align*}
\psi_{\beta,\delta_{-1}}\big(V_aS^mS^{*n}V_a^*\big)&=\zeta(\beta)^{-1}a^{-\beta}\sum_{c\in\nx}(2c)^{-\beta}-\zeta(\beta)^{-1}a^{-\beta}\sum_{c\in\nx\setminus 2\nx}c^{-\beta}\\
&=a^{-\beta}\big(2^{-\beta}-(1-2^{-\beta}\big))\\
&=a^{-\beta}(2.2^{-\beta}-1)=a^{-\beta}(2^{1-\beta}-1).
\end{align*}
To sum up, we have
\[
\psi_{1,\delta_{-1}}\big(V_aS^mS^{*n}V_a^*\big)=
\begin{cases}
a^{-1}&\text{if $m-n$ is even, and}\\
0&\text{if $m-n$ is odd.}
\end{cases}
\]
\end{example}

So we see by example that there are many different KMS$_1$ states of the system $(\TT(\nxxn),\sigma)$. This is of course quite different from what happened at the critical inverse temperature for the left Toeplitz algebra $\TT(\nxnx)$ (see \cite[Theorem~7.1(2)]{LR-advmath}).

\appendix \section{Concluding remarks}\label{rabbit}

\subsection{Finitely many primes}\label{finitep} We now consider a finite set $E$ of prime numbers, and the set $\N_E$ of integers $n\in \N$ whose prime factors all belong to $E$. Then $\N_E$ is an (additive) subgroup of $\N$, and the set $\N_E^\times=\N_E\setminus\{0\}$ is a (multiplicative) subgroup of $\nx$ which is isomorphic via prime factorisation to $\N^E\cong\N^{|E|}$. 

Since the case $d<\infty$ is explicitly allowed in \cite[\S4]{diri}, we can run the argument of \S\ref{sec-KMS-largebeta} for the semidirect product $\N_E^\times\ltimes \N_E$. However, because there are now only finitely many primes in play, the Dirichlet function $\sum_{a\in \N_E^\times} a^{-\beta}$ converges for all $\beta>0$, with sum
\[
\zeta_E(\beta):=\prod_{p\in E} (1-p^{-\beta}).
\]
Hence for all $\beta>0$, there is a bijection $\mu\mapsto \psi_{\mu,\beta}$ of the set $P(\T)$ of probability measures onto the simplex of KMS$_\beta$ states of $(\TT(\N_E^\times\ltimes\N_E),\sigma)$ such that
\begin{equation}\label{comppsi}
\psi_{\mu,\beta}(V_aS^mS^{*n}V_b^*)=\delta_{a,b}\frac{a^{-\beta}}{\zeta_E(\beta)}\sum_{c\in \N_E^\times}c^{-\beta}\int_{\T}z^{c(m-n)}\,d\mu(z).
\end{equation}

\subsection{A crossed-product realisation of the boundary quotient}
The commutation relations in Proposition~\ref{presentaddquot} make the additive boundary quotient $\partial_{\add}\TT(\nxxn)$ look a little like a crossed product of the singly generated $C^*$-algebra $C^*(S)$ by an endomorphic action $\alpha$ of $\nx$ satisfying $\alpha_a(S)=S^a$.  However, it is not one of the familiar ``semigroup crossed products'' studied in \cite{St}, \cite{quasilat}, \cite{bcalg}, \cite{E} or \cite{aHR}, for example; it is more like the ``semicrossed products" introduced recently in \cite{DFK} and \cite{kaka}. Here we make this connection precise. 

We  view $\partial_{\add}\TT(\nxxn)$  as the universal $C^*$-algebra generated by elements satisfying the relations of Proposition~\ref{presentaddquot}. As usual, we use lower case $(s,v_a)$ for the generators of $\partial_{\add}\TT(\nxxn)$ to emphasise that they have a universal property. 

To write $\partial_{\add}\TT(\nxxn)$ as a crossed product, we need to identify the coefficient algebra $C^*(s)$. Relation~(A3) says that $s$ is unitary, and implies that $C^*(s)$ is commutative and isomorphic to $C(\sigma(s))$. Since $s$ is unitary, we have $\sigma(s)\subset \T$, and we claim that it is all of $\T$. To see this, we use the representation $(S,V)$ on $\ell^2(\nx\rtimes \Z)$ of Example~\ref{exDFKcovrep}. An elementary calculation shows that, for each $z\in\T$, the point mass $e_{1,0}$ in $\ell^2(\nx\rtimes \Z)$ is not in the range of $S-z1$, and hence $z\in \sigma(S)$. Since $\sigma(S)\subset\sigma(s)\subset \T$, we deduce that $\sigma(s)$ is also all of $\T$. Thus $(C^*(s),s)$ is canonically isomorphic to $(C(\T),\iota)$, where $\iota$ is the function $\iota:z\mapsto z$.

Since $C(\T)$ is the universal algebra generated by the unitary $\iota$, and each $\iota^a$ is unitary, there are endomorphisms $\{\alpha_a:a\in\nx\}$ of $C(\T)$ characterised by $\alpha_a(\iota)=\iota^a$. It follows by approximation that $\alpha_a$ is given on more general functions $f\in C(\T)$ by $\alpha_a(f)(z)=f(z^a)$. 

We now consider a Nica--Toeplitz--Pimsner algebra $\NT(C(\T),\nx,\alpha)$ which is similar to the ones studied by Kakariadis in \cite{kaka}. Our  $\NT(C(\T),\nx,\alpha)$   is universal for pairs consisting of a unitary operator $S$ on a Hilbert space $H$ and a Nica-covariant representation $V$ of $\nx$ on $H$ such that 
\begin{equation}\label{NCPcov}
\pi(f)V_a=V_a\pi(\alpha_a(f)) \quad\text{for $f\in C(\T)$ and $a\in \nx$.}
\end{equation}

\begin{remark}
Our $\NT(C(\T),\nx,\alpha)$ is similar to, but different from, that of \cite[Definition~2.1]{kaka}.  First, the unitary operator $S\in U(H)$ gives a unital representation $\pi_S$ of the coefficient algebra $A=C(\T)$, and since we can recover $S$ via the formula $S:=\pi_S(\iota)$, these unitaries $S$ are in one-to-one correspondence with the representations $\pi$ described in \cite[Definition~2.1(i)]{kaka}. 
Second, our semigroup $\nx$ is not finitely generated whereas the semigroup $\Z^n_+$ in \cite{kaka} (which we would  denote by $\N^n$) is finitely generated.   If we fix $n$ primes $\{p_i:1\leq i\leq n\}$, then $V\mapsto \{V_{p_i}:1\leq i\leq n\}$  gives a generating family $\{V_i:1\leq i \leq n\}$ as in \cite[Definition~2.1(ii)]{kaka}. Restricting our covariance relation~\eqref{NCPcov} to this family gives the relation in \cite[Definition~2.1(iii)]{kaka}. 
\end{remark}

\begin{prop}\label{=kakacp}
Let $\alpha:\nx\to \End C(\T)$ be the action characterised by $\alpha_a(\iota)=\iota^a$. Then the additive boundary quotient $\partial_{\add}\TT(\nxxn)$ is the Toeplitz--Nica--Pimsner algebra $\NT(C(\T),\nx,\alpha)$ defined at (\ref{NCPcov}) and the boundary quotient $\partial\TT(\nxxn)$ is the Cuntz--Nica--Pimsner algebra $\NO(C(\T),\nx,\alpha)$ analogous to the one from \cite[\S2]{kaka}. 
\end{prop}

\begin{proof}
Proposition~\ref{presentaddquot} says that $\partial_{\add}$ and $\NT(C(\T),\nx,\alpha)$ have equivalent presentations, and hence are canonically isomorphic algebras.

For the second assertion, we need to check the extra relations imposed in \cite[\S2.3]{kaka}. View $\nx$ as $\N^{\PP}$. Since each $z\mapsto z^a$ is surjective, each endomorphism $\alpha_p$ is injective, and the ideals $I_p$ in \cite{kaka} are all $A=C(\T)$. So the extra relations are 
\[
\prod_{p\in S}(1-V_pV_p^*)=0\quad\text{for all finite subsets $S$ of $\PP$.}
\]
These are equivalent to $1-V_pV_p^*=0$ for all $p\in \PP$, and hence say merely that each $V_a$ is unitary. Thus the isomorphism $\pi_S\times V$ induces an isomorphism of the Crisp--Laca boundary quotient $\partial TT(\nxxn)$ onto the quotient $\NO(C(\T),\nx,\alpha)$ of $\NT(C(\T),\nx,\alpha)$.
\end{proof}

\begin{remark}
It follows from Proposition~\ref{=kakacp} and Proposition~\ref{propCL=cp} that the Cuntz--Nica--Pimsner algebra $\NO(C(\T),\nx,\alpha)$ is isomorphic to the group algebra $C^*(\qxxq)$. Indeed, the algebra $\NO(C(\T),\nx,\alpha)$ is generated by a unitary element $S$ and a unitary representation $U$, and the proof of \S\ref{strucbdaryq} shows how such $S$ and $U$ combine to give a unitary representation $W$ of $\Q_+$ (and hence also of $\Q$) alongside the representation $U$ of ${\qx}$ coming from $V$. The pair $(U,W)$ now give a unitary representation of $\qxxq$, and hence a homomorphism of $C^*(\qxxq)$ into $\NO(C(\T),\nx,\alpha)$ which is the required isomorphism. 
\end{remark}

\subsection{Relations with work of Kakariadis on KMS states} In \cite[\S4]{kaka}, Kakariadis has computed the KMS states on crossed products of the the form $\NT(A,\N^n,\alpha)$ for endomorphic actions $\alpha$ of the additive semigroup $\N^n$ (denoted there by $\Z^n_+$) on unital $C^*$-algebras $A$. Our results overlap with his, and especially in the case of finitely many primes, as discussed in \S\ref{finitep}, and $A=C(\T)$. 

To see this, we suppose again that $E$ is a finite set of primes, and use the notation of \S\ref{finitep}. Then the map $M:\N_E\to \N_E^\times$ defined by $M(x)=\prod_{p\in E}p^{x_p}$ is an isomorphism. We define $\tau:\N_E\to \End C(\T)$ in terms of the action $\alpha:\N_E^\times\to \End C(\T)$ of Proposition~\ref{=kakacp} by $\tau=\alpha\circ M$. We then define $\lambda\in [0,\infty)^E$ by $\lambda_p=\ln p$, and use it to define an action $\sigma^\lambda$ of $\R$ on $\NT(C(\T),\N_E,\tau)$, as in \cite[\S4.1]{kaka}. (We have simplified his notation a little by noting that $\lambda_p\beta>0$ for all $p$, and hence there is no need to put hypotheses on the vector $\underline{\beta}$.)
The simplex of tracial states on $C(\T)$ is the simplex of $P(\T)$ of probability measures on $\T$. Hence \cite[Theorem~4.4]{kaka} says that, for every $\beta>0$, there is an affine homeomorphism $\mu\mapsto\phi_{\mu,\beta}$ of $P(\T)$ onto the simplex of KMS$_\beta$ states of $\NT(C(\T),\N_E,\tau)$ such that
\[
\phi_{\mu,\beta}(V_xfV_y^*)=\delta_{x,y}e^{-\langle x,\lambda\rangle\beta}\prod_{p\in E}(1-e^{-\beta\lambda_p})\Big(\sum_{w\in \N_E} e^{-\langle w,\lambda\rangle\beta}\int_{\T}\rho_w(f)\,d\mu\Big)
\]
for all $x,y\in\N_E$ and $f\in C(\T)$.

To transfer this into our notation, we note that
\[
\langle x,\lambda\rangle=\sum_{p\in E} x_p\lambda_p=\sum_{p\in E} x_p\ln p,
\]
and hence
\[
e^{-\langle x,\lambda\rangle\beta}=\prod_{p\in E}e^{-(x_p\ln p)\beta}=\prod_{p\in E}p^{-x_p\beta}=M(x)^{-\beta}.
\]
Similarly, we have $e^{-\beta\lambda_p}=p^{-\beta}$. Next we take $f(z)=z^m$ and note that
\[
\tau_w(f)(z)=\alpha_{M(w)}(f)(z)=f(z^{M(w)})=z^{M(w)m}.
\]
Thus we have
\[
\phi_{\mu,\beta}(V_xz^mV_y^*)=\delta_{x,y} M(x)^{-\beta}\prod_{p\in E}(1-p^{_\beta})
\Big(\sum_{w\in \N_E} M(w)^{-\beta}\int_{\T}z^{M(w)m}\,d\mu(z)\Big).
\]
When we pull this back using the isomorphism $M:\N_E\to N_E^\times$ to a state on $\NT(C(\T),\N_E^\times,\alpha)$ (effectively replacing $x,y$ by $a=M(x),b=M(y)$), we recover the state $\psi_{\mu,\beta}$ of \S\ref{comppsi}.

\subsection{LCM semigroups} Afsar, Brownlowe, Larsen and Stammeier have recently proved a very general result \cite[Theorem~4.3]{ABLS} about the KMS states on the $C^*$-algebras of right LCM semigroups which unifies the main results of \cite{LR-advmath}, \cite{BaHLR}, \cite{LRR}, \cite{LRRW} and \cite{CaHR}. So it is natural to ask whether their result applies also to the semigroup $\nxxn$. The answer is that it does not seem to. To see why, we need to view $\nxxn$ as a right LCM semigroup, and see that it is not ``admissible'' in the sense of \cite[Definition~3.1]{ABLS}.

The semigroup $S:=\nxxn$ is certainly a right LCM semigroup, because it is the positive cone in the quasi-lattice ordered group $(\qxxq,\nxxn)$. Indeed, in the original development \cite[\S3--4]{BLS} of the theory of right LCM semigroup algebras, the authors  used the paper \cite{quasilat} on quasi-lattice ordered groups as a template. So we can try to apply the result in \cite{ABLS} to $S=\nxxn$. (As a trivial point of notation,  we continue to write $s\vee t<\infty$ to mean $sS\cap tS\not=\emptyset$.) Since $S$ is lattice ordered, the ``core subsemigroup''
\[
S_c:=\big\{s\in S:s\vee t<\infty\text{ for all $t\in S$}\big\}
\]
is all of $S$. Thus there are no elements of $S\backslash S_c$, hence trivially no ``core irreducible'' elements. In the notation of \cite{ABLS}, we have $S_{ci}^1=\{1\}$ and $S_{ci}^1S_c=S$, which is (A1) in  \cite[Definition~3.1]{ABLS}. The second condition (A2) holds vacuously. But (A3), which postulates the existence of a homomorphism $N:S\to \nx$ with certain properties, seems to be a problem.

Suppose that $N:\nxxn\to \nx$ is a homomorphism. Then we take $a\in \nx$ with $a>1$, and compute $N(a,a)$ two ways. On the one hand, we have
\[
N(a,a)=N(\big((1,1)(a,0)\big)=N(1,1)N(a,0), 
\]
and on the other we have
\[
N(a,a)=N\big((a,0),(1,a)\big)=N(a,0)N\big((1,1)^a\big)=N(a,0)N(1,1)^a.
\]
Since $a>1$, we deduce that $N(1,1)=1$, and
\[
N(a,m)=N\big((a,0)(1,m)\big)=N(a,0)N(1,1)^m=N(a,0)\quad\text{for all $(a,m)\in \nxxn$.}
\]
Thus for $c\in \nx$ belonging to the range of $N$, the set $N^{-1}(c)$ contains all $(a,m)$ with $N(a,0)=c$, and is in particular infinite. But here, since $S_c=S$, the equivalence relation $\sim$ used in \cite{ABLS} satisfies $s\sim t\Longleftrightarrow s\vee t<\infty$, and hence we have $s\sim t$ for all $s,t$ in our lattice-ordered semigroup $S$. Thus $|N^{-1}(c)/\!\sim|=1$, and $N$ cannot satisfy condition (A3)(a) in  \cite[Definition~3.1]{ABLS}. So Theorem~4.3 of \cite{ABLS} does not apply to the semigroup $\nxxn$.

\end{document}